\input amstex
\documentstyle{amsppt}
%
\catcode`@=11
\redefine\output@{%
  \def\break{\penalty-\@M}\let\par\endgraf
  \ifodd\pageno\global\hoffset=105pt\else\global\hoffset=8pt\fi  
  \shipout\vbox{%
    \ifplain@
      \let\makeheadline\relax \let\makefootline\relax
    \else
      \iffirstpage@ \global\firstpage@false
        \let\rightheadline\frheadline
        \let\leftheadline\flheadline
      \else
        \ifrunheads@ 
        \else \let\makeheadline\relax
        \fi
      \fi
    \fi
    \makeheadline \pagebody \makefootline}%
  \advancepageno \ifnum\outputpenalty>-\@MM\else\dosupereject\fi
}
\catcode`\@=\active
\nopagenumbers
\def\negskp{\hskip -2pt}
\def\Img{\operatorname{Im}}
\def\Ker{\operatorname{Ker}}
\def\const{\operatorname{const}}
\def\compos{\,\raise 1pt\hbox{$\sssize\circ$} \,}
\def\idop{\operatorname{\bold{id}}}
\def\msum#1{\operatornamewithlimits{\sum^#1\!{\ssize\ldots}\!\sum^#1}}
\accentedsymbol\tx{\tilde x}
\accentedsymbol\hbh{\kern-3pt\hat{\kern 3pt\bold h}}
\accentedsymbol\tbb{\kern-3pt\tilde{\kern 3pt\bold b}}
\accentedsymbol\vnabla{\nabla\kern -7pt\raise 5pt\vbox{\hrule width 5.5pt}
\kern 1.7pt}
\def\blue#1{#1}
\catcode`#=11\def\diez{#}\catcode`#=6
\catcode`_=11\def\podcherkivanie{_}\catcode`_=8
\catcode`~=11\catcode`~=\active
\def\mycite#1{\cite{\blue{#1}}\immediate\special{ps:
     ShrHPSdict begin /ShrBORDERthickness 0 def}}
\def\myciterange#1#2{\cite{\blue{#2}}\immediate\special{ps:
     ShrHPSdict begin /ShrBORDERthickness 0 def}}
\def\mytag#1{%
    \tag#1}
\def\mythetag#1{\thetag{\blue{#1}}\immediate\special{ps:
     ShrHPSdict begin /ShrBORDERthickness 0 def}}
\def\myrefno#1{\no#1}
\def\myhref#1#2{\blue{#2}\immediate\special{ps:
     ShrHPSdict begin /ShrBORDERthickness 0 def}}
\def\myEarXivlink{\myhref{http://arXiv.org}{http:/\negskp/arXiv.org}}
\def\myGeoCities{\myhref{http://www.geocities.com}{GeoCities}}
\def\mytheorem#1{\csname proclaim\endcsname{Theorem #1}}
\def\mythetheorem#1{\blue{#1}\immediate\special{ps:
     ShrHPSdict begin /ShrBORDERthickness 0 def}}
\def\mylemma#1{\csname proclaim\endcsname{Lemma #1}}
\def\mythelemma#1{\blue{#1}\immediate\special{ps:
     ShrHPSdict begin /ShrBORDERthickness 0 def}}
\def\mydefinition#1{\definition{Definition #1}}
\def\mythedefinition#1{\blue{#1}\immediate\special{ps:
     ShrHPSdict begin /ShrBORDERthickness 0 def}}
\def\myanchortext#1#2{\blue{#2}}
\def\mytheanchortext#1#2{\blue{#2}\immediate\special{ps:
     ShrHPSdict begin /ShrBORDERthickness 0 def}}
\font\tencyr=wncyr10
\font\eightcyr=wncyr8
\pagewidth{360pt}
\pageheight{606pt}
\topmatter
\title
Tensor functions of tensors\\
and the concept of extended tensor fields.
\endtitle
\author
Ruslan Sharipov
\endauthor
\address Rabochaya street 5, 450003 Ufa, Russia
\endaddress
\email \vtop to 30pt{\hsize=280pt\noindent
\myhref{mailto:R\podcherkivanie Sharipov\@ic.bashedu.ru}
{R\_\hskip 1pt Sharipov\@ic.bashedu.ru}\newline
\myhref{mailto:r-sharipov\@mail.ru}
{r-sharipov\@mail.ru}\newline
\myhref{mailto:ra\podcherkivanie sharipov\@lycos.com}{ra\_\hskip 1pt
sharipov\@lycos.com}\vss}
\endemail
\urladdr
\vtop to 20pt{\hsize=280pt\noindent
\myhref{http://www.geocities.com/r-sharipov}
{http:/\negskp/www.geocities.com/r-sharipov}\newline
\myhref{http://www.freetextbooks.boom.ru/index.html}
{http:/\negskp/www.freetextbooks.boom.ru/index.html}\vss}
\endurladdr
\abstract
    Tensor fields depending on other tensor fields are considered.
The concept of extended tensor fields is introduced and the theory
of differentiation for such fields is developed.
\endabstract
\endtopmatter
\loadbold
\TagsOnRight
\document

\rightheadtext{Tensor functions of tensors \dots}
\head
1. Tensors and tensor fields on manifolds.
\endhead
    Let $M$ be some $n$-dimensional smooth real manifold. Then each
point $p\in M$ has some neighborhood $U$ bijectively mapped onto an
open set $V$ in the $n$-dimensional space $\Bbb R^n$. This means that
any point $p\in U$ is associated with some unique vector in $V$ with the 
coordinates $x^1(p),\,\ldots,\,x^n(p)$. The set $V\subset\Bbb R^n$ is
called a {\it local map\/} or a {\it local chart\/} of the manifold $M$,
while the numbers $x^1(p),\,\ldots,\,x^n(p)$ are called the {\it 
coordinates\/} of the point $p$ in the local chart $V$. In a not too
formal terminology the set $U\subset M$ is also called a {\it local 
map\/} or a {\it local chart}.\par
     The whole manifold is covered by local charts. If two local 
charts $U$ and $\tilde U$ do overlap, i\.\,e\. if \ $U\cap \tilde
U\neq\varnothing$, then the so-called {\it transition functions\/} arise:
$$
\xalignat 2
&\hskip -2em
\cases
\tx^1=\tx^1(x^1,\,\ldots,x^n),\\
. \ . \ . \ .\ . \ . \ . \ . \ . \ . \ 
. \ . \ . \ . \ .\\
\tx^n=\tx^n(x^1,\,\ldots,x^n),
\endcases
&&\cases
x^1=x^1(\tx^1,\,\ldots,\tx^n),\\
. \ . \ . \ .\ . \ . \ . \ . \ . \ . \ 
. \ . \ . \ . \ .\\
x^n=x^n(\tx^1,\,\ldots,\tx^n).
\endcases
\mytag{1.1}
\endxalignat
$$
They relate the local coordinates of a point $p\in U\cap \tilde U$ in
two charts. In the case of smooth manifolds the transition functions
\mythetag{1.1} for all pairs of overlapping maps are smooth functions. 
The partial derivatives of \mythetag{1.1} form the transition matrices:
$$
\xalignat 2
&\hskip -2em
S^i_j=\frac{\partial x^i}{\partial\tx^j},
&&T^i_j=\frac{\partial\tx^i}{\partial x^j}.
\mytag{1.2}
\endxalignat
$$
They are inverse to each other: $S=T^{-1}$. By tradition $S$ is called 
the {\it direct transition matrix}, while $T$ is called the {\it inverse
transition matrix}.\par
     All what was said just above is a standard definition of a smooth 
real manifold. We give it here in order to make this paper understandable
not only to professional mathematicians, but to physicists and to students
majoring in physics and engineering. With the same purpose in mind, below 
we shall combine the coordinate and coordinate-free approaches.\par
    Continuing our introductory section, let's consider the following
differential\linebreak operators associated with two local coordinate
systems in $M$:
$$
\xalignat 2
&\hskip -2em
\bold E_i=\frac{\partial}{\partial x^i},
&&\tilde\bold E_i=\frac{\partial}{\partial\tx^i}.
\mytag{1.3}
\endxalignat
$$
From \mythetag{1.1} and \mythetag{1.2} it is easy to derive that
the differential operators $\bold E_i$ and $\tilde\bold E_j$ are 
related to each other through the following equalities:
$$
\xalignat 2
&\hskip -2em
\tilde\bold E_j=\sum^n_{i=1}S^i_j\ \bold E_i,
&&\bold E_i=\sum^n_{j=1}T^j_i\ \tilde\bold E_j.
\mytag{1.4}
\endxalignat
$$\par
    The three-dimensional Euclidean space $\Bbb E$, which we observe 
in our everyday life, is an example of a smooth manifold. In this case 
$x^1,\,\ldots,\,x^n$ and $\tx^1,\,\ldots,\,\tx^n$ are interpreted as
Cartesian and/or curvilinear coordinates in $\Bbb E$. The differential
operators \mythetag{1.3} can be associated with the frame vectors of
moving frames of two curvilinear coordinate systems (see \mycite{1}
and \mycite{2}) because \mythetag{1.4} coincide with the corresponding
relationships for the frame vectors.\par
    A smooth surface in the space $\Bbb E$ is another example of a
smooth manifold. In this case the differential operators \mythetag{1.3}
can be associated with the tangent vectors forming a basis in the tangent
plane to that surface (see \mycite{2}). This is the reason why in the
case of an arbitrary smooth manifold $M$ the operators \mythetag{1.3} are 
called {\it tangent vectors}. If some point $p\in U$ is fixed, then the
{\it tangent space\/} $T_p(M)$ is defined as the span of the vectors
\mythetag{1.3} at that point:
$$
\hskip -2em
T_p(M)=\bigl<\bold E_1,\,\ldots,\,\bold E_n\bigr>
=\bigl<\tilde\bold E_1,\,\ldots,\,\tilde\bold E_n\bigr>.
\mytag{1.5}
$$
The tangent spaces $T_p(M)$ and $T_q(M)$ of different points $p\neq q$
are understood as two different $n$-dimensional vector spaces\footnotemark.
\footnotetext{ \ Informally, one can imagine a manifold $M$ as a 
\myanchortext{1}{cat} with the hairs $T_p(M)$ growing from each point $p$
on its skin.}
\adjustfootnotemark{-1}
\mydefinition{1.1} A vector field $\bold X$ is a vector-valued function
that maps each point $p\in M$ to some vector $\bold X(p)\in T_p(M)$.
\enddefinition
    This is an invariant (coordinate-free) definition of a vector field.
Due to \mythetag{1.5} one can expand the vector $\bold X(p)$ in two different bases:
$$
\xalignat 2
&\hskip -2em
\bold X(p)=\sum^n_{i=1}X^i\ \bold E_i,
&&\bold X(p)=\sum^n_{j=1}\tilde X^j\ \tilde\bold E_i.
\mytag{1.6}
\endxalignat
$$
Here $X^i=X^i(p)=X^i(x^1,\ldots,x^n)$ and $\tilde X^i=\tilde X^i(p)
=\tilde X^i(\tx^1,\ldots,\tx^n)$ are the components of the vector 
$\bold X(p)$ in two different local charts. Due to \mythetag{1.4} and
\mythetag{1.6} they are related to each other according to the formulas:
$$
\xalignat 2
&\hskip -2em
X^i=\sum^n_{j=1}S^i_j\ \tilde X^j,
&&\tilde X^j=\sum^n_{i=1}T^j_i\ X^i.
\mytag{1.7}
\endxalignat
$$
The formulas \mythetag{1.7} form the base for coordinate definition of a
vector field (see \mycite{1} and \mycite{2}). Here this definition is
formulated as follows.
\mydefinition{1.2} A vector field $\bold X$ is a geometric object 
in each local chart represented by its components $X^i=X^i(x^1,\ldots,
x^n)$ and such that under a change of a local chart its components are
transformed according to the formulas \mythetag{1.7}.
\enddefinition
    Let $T^*_p(M)$ be the dual space for the tangent space
$T_p(M)$ (see \mycite{3} for the definition of a dual space). Then consider
the following tensor product\footnotemark:
\adjustfootnotemark{-1}
$$
\hskip -2em
T^r_s(p,M)=\overbrace{T_p(M)\otimes\ldots\otimes T_p(M)}^{\text{$r$
times}}\otimes\underbrace{T^*_p(M)\otimes\ldots\otimes T^*_p(M)}_{\text{$s$
times}}.
\mytag{1.8}
$$
The tensor product \mythetag{1.8} is also a vector space associated with
the point $p$. So, each point of a smooth manifold carries a great many
mathematical constructs\footnotemark\ including but not limited to those considered in this paper. Some of these constructs correspond to real
physical fields, others are waiting their time to be associated with something in the nature.\par
\footnotetext{ \ One can imagine $T^*_p(M)$ and $T^r_s(p,M)$ as other hairs
on the skin of that our \mytheanchortext{1}{cat} growing from the same point
$p$ as $T_p(M)$. However, I don't know if some real animal can have a bunch
of hairs on the same root.}
\mydefinition{1.3} A tensor field $\bold X$ of the type $(r,s)$ is a
tensor-valued function that maps each point $p\in M$ to some tensor 
$\bold X(p)\in T^r_s(p,M)$.
\enddefinition
    In order to represent a tensor field in a local map we need to have
some basis in the space \mythetag{1.8}. The differentials $dx^1,\,\ldots,
\,dx^n$ form a basis in the dual space $T^*_p(M)$. This basis is dual to
the basis of tangent vectors $\bold E_1,\,\ldots,\,\bold E_n$ in $T_p(M)$:
$$
\hskip -2em
\bigl<dx^i\,|\,\bold E_j\bigr>=\delta^i_j=\cases 1 &\text{for \ }i=j,\\
0&\text{for \ }i\neq j.
\endcases
\mytag{1.9}
$$
By angular brackets in \mythetag{1.9} we denote the scalar product
of a vector and a covector (see definition in Chapter~\uppercase
\expandafter{\romannumeral 3} of \mycite{3}). Now let's denote
$$
\hskip -2em
\bold E^{j_1\ldots\,j_s}_{i_1\ldots\,i_r}=\bold E_{i_1}
\otimes\ldots\otimes\bold E_{i_r}\otimes dx^{j_1}\otimes
\ldots\otimes dx^{j_s}.
\mytag{1.10}
$$
The local\footnotemark\ tensor fields \mythetag{1.10} form a basis in the
space \mythetag{1.8} at all points $p\in U$. Therefore, any tensor field
$\bold X$ of the type $(r,s)$ admits the expansion
$$
\hskip -2em
\bold X=\sum^n_{i_1=1}\ldots\sum^n_{i_r=1}\sum^n_{j_1=1}\ldots
\sum^n_{j_s=1}X^{i_1\ldots\,i_r}_{j_1\ldots\,j_s}\ 
\bold E^{j_1\ldots\,j_s}_{i_1\ldots\,i_r}.
\mytag{1.11}
$$
\footnotetext{ \ The tensor fields \mythetag{1.10} are defined only within
the local chart $U$.}\adjustfootnotemark{-2}The coefficients $X^{i_1\ldots
\,i_r}_{j_1\ldots\,j_s}=X^{i_1\ldots\,i_r}_{j_1\ldots\,j_s}(x^1,\,\ldots,
x^n)$ in the expansion \mythetag{1.11} are called the {\it components\/} 
of the tensor field $\bold X$ in the local chart $U$. Under a change of a
local chart they are transformed as follows:
$$
\pagebreak
\hskip -2em
X^{i_1\ldots\,i_r}_{j_1\ldots\,j_s}=\msum{n}\Sb h_1,\,\ldots,\,h_r\\
k_1,\,\ldots,\,k_s\endSb S^{i_1}_{h_1}\ldots\,S^{i_r}_{h_r}\ 
T^{k_1}_{j_1}\ldots\,T^{k_s}_{j_s}\ \tilde X^{h_1\ldots\,h_r}_{k_1
\ldots\,k_s},
\mytag{1.12}
$$
The transformation rule can be inverted. The inverse transformation is
written as
$$
\hskip -2em
\tilde X^{i_1\ldots\,i_r}_{j_1\ldots\,j_s}=\msum{n}\Sb h_1,\,\ldots,\,h_r\\
k_1,\,\ldots,\,k_s\endSb T^{i_1}_{h_1}\ldots\,T^{i_r}_{h_r}\ S^{k_1}_{j_1}
\ldots\,S^{k_s}_{j_s}\ X^{h_1\ldots\,h_r}_{k_1\ldots\,k_s}.
\mytag{1.13}
$$
An alternative definition of a tensor field is based on the transformation
rules \mythetag{1.12} and \mythetag{1.13} for its components.
\mydefinition{1.4} A tensor field $\bold X$ of the type $(r,s)$ is a
geometric object in each local chart represented by its components 
$X^{i_1\ldots\,i_r}_{j_1\ldots\,j_s} =X^{i_1\ldots\,i_r}_{j_1\ldots
\,j_s}(x^1,\,\ldots,x^n)$ and such that under a change of a local chart
its components obey the transformation\linebreak rules \mythetag{1.12} 
and \mythetag{1.13}.
\enddefinition
\head
2. Tangent bundle, cotangent bundle,\\
and other tensor bundles.
\endhead
     Let $p$ be some point of a smooth real manifold $M$ and let $\bold v
\in T_p(M)$ be some tangent vector at the point $p$. The set of all pairs
$q=(p,\bold v)$ forms another smooth real manifold. It is called the
{\it tangent bundle\footnotemark\/} of $M$ and denoted by $TM$. The map
$\pi\!:TM\to M$ that takes a point $q=(p,\bold v)$ of $TM$ to the point
$p\in M$ is called the {\it canonical projection\/} of the tangent bundle
$TM$ onto the base manifold $M$. If a local chart $U$ on $M$ is given, a
point $q=(p,\bold v)$ of $TM$ is represented by $2n$ variables
$$
\hskip -2em
x^1,\,\ldots,\,x^n,\,v^1,\,\ldots,\,v^n,
\mytag{2.1}
$$
\footnotetext{ \ Being more strict, $TM$ is called the {\it total 
space\/} of a tangent bundle, while a {\it tangent bundle\/} itself 
is a whole construct including a total space, a base, and a projection 
map $\pi\!:TM\to M$. However, in this paper we use more loose 
terminology.}\adjustfootnotemark{-1}where $x^1,\,\ldots,\,x^n$ are the
local coordinates of the point $p=\pi(q)$ and $v^1,\,\ldots,\,v^n$ are 
the components of the tangent vector $\bold v$:
$$
\hskip -2em
\bold v=v^1\ \bold E_1+\ldots+v^n\ \bold E_n.
\mytag{2.2}
$$
Hence, we have $\dim(TM)=2\,\dim(M)$. Tangent bundles naturally arise
in considering Newtonian dynamical systems with holonomic constraints 
(see theses \mycite{4}, \mycite{5}, and the series of papers  
\myciterange{6}{6--23}). In mechanics a base manifold $M$ is called
a {\it configuration space}, $TM$ is called a {\it phase space}, and
\mythetag{2.2} is interpreted as the {\it velocity vector} of a point
moving within $M$.\par
    Let $\bold p\in T^*_p(M)$ be some covector at the point $p\in M$.
The set of all pairs $q=(p,\bold p)$ forms a smooth real manifold
which is called the {\it cotangent bundle\/} of $M$. It is denoted
$T^*\!M$. In the case of cotangent bundle $T^*\!M$ we also have the 
{\it canonical projection\/} $\pi\!:T^*\!M\to M$ that takes a point
$q=(p,\bold p)$ to the point $p$ of the base manifold $M$. In a local 
chart a point $q=(p,\bold p)$ is represented by $2n$ variables
$$
\hskip -2em
x^1,\,\ldots,\,x^n,\,p_1,\,\ldots,\,p_n,
\mytag{2.3}
$$
where $x^1,\,\ldots,\,x^n$ are the local coordinates of the point
$p=\pi(q)$ and $p_1,\,\ldots,\,p_n$ are the components of the 
covector $\bold p\in T^*_p(M)$:
$$
\hskip -2em
\bold p=p_1\ dx^1+\ldots+p_n\ dx^n.
\mytag{2.4}
$$
For the dimension of a cotangent bundle we have $\dim(T^*\!M)=2\,
\dim(M)$. Cotangent bundles naturally arise when one passes from 
Lagrangian dynamical systems to the corresponding Hamiltonian 
dynamical systems. In this case $T^*\!M$ is interpreted as the 
$\bold p$-representation of a phase space $TM$, while $TM$ is 
understood as the $\bold v$-representation of $T^*\!M$ (see papers
\myciterange{24}{24--30}). In Hamiltonian dynamics the covector \mythetag{2.4} is called the {\it momentum covector}.\par
     Tensor bundles are defined by analogy to $TM$ and $T^*\!M$.
Let $p$ be some point of the base manifold $M$ and let $\bold T$
be some tensor of the type $(r,s)$ at this point. The set of all
pairs $q=(p,\bold T)$ forms a smooth real manifold $T^r_sM$ of 
the dimension $\dim(T^r_sM)=n+n^{r+s}$. It is called the 
{\it tensor bundle\/} of the type $(r,s)$ over the base $M$.
The map $\pi\!:T^r_sM\to M$ that takes a point $q=(p,\bold T)$ of 
$T^r_sM$ to the point $p\in M$ is called the {\it canonical 
projection\/} of the tensor bundle $T^r_sM$ onto the base manifold 
$M$. Like in \mythetag{2.1} and \mythetag{2.3}, we can specify the
set of variables associated with a point $q=(p,\bold T)$ of the
tensor bundle $T^r_sM$ and with a local chart $U\subset M$:
$$
\hskip -2em
x^1,\,\ldots,\,x^n,\,T^{1\kern 1pt\ldots\,1}_{1\kern 1pt\ldots\,1},
\,\ldots,\,T^{n\,\ldots\,n}_{n\,\ldots\,n}.
\mytag{2.5}
$$
The number of variables in \mythetag{2.5} determines the dimension
$\dim(T^r_sM)=n+n^{r+s}$ of the tensor bundle $T^r_sM$. Note that if
we consider a tensor field $\bold T$, its components
$T^{i_1\ldots\,i_r}_{j_1\ldots\,j_s}=T^{i_1\ldots\,i_r}_{j_1\ldots
\,j_s}(x^1,\,\ldots,x^n)$ are functions of $x^1,\,\ldots,\,x^n$, while
$T^{i_1\ldots\,i_r}_{j_1\ldots\,j_s}$ in \mythetag{2.5} are independent
variables. For another chart $\tilde U\subset M$ we have the other
set of variables
$$
\hskip -2em
\tx^1,\,\ldots,\,\tx^n,\,\tilde T^{1\kern 1pt\ldots\,1}_{1\kern
1pt\ldots\,1},\,\ldots,\,\tilde T^{n\,\ldots\,n}_{n\,\ldots\,n}.
\mytag{2.6}
$$
If the charts $U$ and $\tilde U$ are overlapping, then we can write
the transition functions
$$
\hskip -2em
\cases
\tx^1=\tx^1(x^1,\,\ldots,x^n),\\
. \ . \ . \ .\ . \ . \ . \ . \ . \ . \ 
. \ . \ . \ . \ .\\
\tx^n=\tx^n(x^1,\,\ldots,x^n),\\
\tilde T^{i_1\ldots\,i_r}_{j_1\ldots\,j_s}=\dsize\msum{n}\Sb h_1,\,\ldots,\,h_r\\
k_1,\,\ldots,\,k_s\endSb T^{i_1}_{h_1}\ldots\,T^{i_r}_{h_r}\ S^{k_1}_{j_1}
\ldots\,S^{k_s}_{j_s}\ T^{h_1\ldots\,h_r}_{k_1\ldots\,k_s},
\endcases
\mytag{2.7}
$$
where the transition matrices $T$ and $S$ are the same as in
\mythetag{1.13} --- they are determined by \mythetag{1.2}. Here
are the inverse transition functions relating \mythetag{2.5} and 
\mythetag{2.6}:
$$
\hskip -2em
\cases
x^1=x^1(\tx^1,\,\ldots,\tx^n),\\
. \ . \ . \ .\ . \ . \ . \ . \ . \ . \ 
. \ . \ . \ . \ .\\
x^n=x^n(\tx^1,\,\ldots,\tx^n),\\
T^{i_1\ldots\,i_r}_{j_1\ldots\,j_s}=\dsize\msum{n}\Sb h_1,\,\ldots,\,h_r\\
k_1,\,\ldots,\,k_s\endSb S^{i_1}_{h_1}\ldots\,S^{i_r}_{h_r}\ T^{k_1}_{j_1}
\ldots\,T^{k_s}_{j_s}\ \tilde T^{h_1\ldots\,h_r}_{k_1\ldots\,k_s}.
\endcases
\mytag{2.8}
$$
The transformations \mythetag{2.7} and \mythetag{2.8} play the same role
for the tensor bundle $T^r_sM$ as the transformations \mythetag{1.1} for
the base manifold $M$. Note that they are linear with respect to
$T^{h_1\ldots\,h_r}_{k_1\ldots\,k_s}$ and $\tilde T^{h_1\ldots\,h_r}_{k_1
\ldots\,k_s}$, but they are nonlinear with respect to the base variables
$x^1,\,\ldots,\,x^n$ and $\tx^1,\,\ldots,\,\tx^n$.
\head
3. Composite tensor bundles.
\endhead
      Suppose that we have several tensor bundles over the same base
manifold $M$. Let's denote them $T^{r_1}_{s_1}\!M,\,\ldots,\,T^{r_Q}_{s_Q}
\kern -1pt M$ and consider their direct sum over $M$:
$$
\hskip -2em
T^{r_1\ldots\,r_Q}_{s_1\ldots\,s_Q}\!M=T^{r_1}_{s_1}\!M\oplus\ldots\oplus
T^{r_Q}_{s_Q}\kern -1pt M
\mytag{3.1}
$$
(see the definition of a direct sum in \mycite{31}). We shall call 
\mythetag{3.1} the {\it composite tensor bundle}. A point $q$ of the
composite tensor bundle \mythetag{3.1} is a list
$$
\hskip -2em
q=(p,\,\bold T[1],\,\ldots,\,\bold T[Q]),
\mytag{3.2}
$$
where $p$ is a point of the base $M$ and $\bold T[1],\,\ldots,\,\bold T[Q]$
are some tensors of the types $(r_1,s_1),\,\ldots,\,(r_Q,s_Q)$ at the point 
$p$. The {\it canonical projection\/} 
$$
\pi\!:\,T^{r_1\ldots\,r_Q}_{s_1\ldots\,s_Q}\!M\to M
$$ is defined as a map that takes a point $q$ of the form \mythetag{3.2} 
to the point $p\in M$. If a local chart $U$ in the base manifold $M$ is
given, a point  $q\in T^{r_1\ldots\,r_Q}_{s_1\ldots\,s_Q}\!M$ such that
$\pi(q)\in U$ is represented by the following set of variables:
$$
x^1,\,\ldots,\,x^n,\,T^{1\kern 1pt\ldots\,1}_{1\kern 1pt\ldots\,1}[1],
\,\ldots,\,T^{n\,\ldots\,n}_{n\,\ldots\,n}[1],\,\ldots,\,
T^{1\kern 1pt\ldots\,1}_{1\kern 1pt\ldots\,1}[Q],
\,\ldots,\,T^{n\,\ldots\,n}_{n\,\ldots\,n}[Q].\qquad
\mytag{3.3}
$$
Taking another chart $\tilde U$, we get the other set of variables 
$$
\tx^1,\,\ldots,\,\tx^n,\,\tilde T^{1\kern 1pt\ldots\,1}_{1\kern
1pt\ldots\,1}[1],\,\ldots,\,\tilde T^{n\,\ldots\,n}_{n\,\ldots
\,n}[1],\,\ldots,\,\tilde T^{1\kern 1pt\ldots\,1}_{1\kern 1pt
\ldots\,1}[Q],\,\ldots,\,\tilde T^{n\,\ldots\,n}_{n\,\ldots\,n}[Q].
\qquad
\mytag{3.4}
$$
If these charts are overlapping, then we have a system of transition functions
$$
\hskip -2em
\cases
\tx^1=\tx^1(x^1,\,\ldots,x^n),\\
. \ . \ . \ .\ . \ . \ . \ . \ . \ . \ 
. \ . \ . \ . \ .\\
\tx^n=\tx^n(x^1,\,\ldots,x^n),\\
\tilde T^{i_1\ldots\,i_r}_{j_1\ldots\,j_s}[P]=\dsize\msum{n}\Sb
h_1,\,\ldots,\,h_r\\
k_1,\,\ldots,\,k_s\endSb T^{i_1}_{h_1}\ldots\,T^{i_r}_{h_r}
\ S^{k_1}_{j_1}\ldots\,S^{k_s}_{j_s}\ T^{h_1\ldots\,h_r}_{k_1
\ldots\,k_s}[P],
\endcases
\mytag{3.5}
$$
where $r=r_P$, $s=s_P$, and the integer number $P$ runs from $1$ to $Q$.
Similarly, we can write the system of the inverse transition functions
for \mythetag{3.5}:
$$
\hskip -2em
\cases
x^1=x^1(\tx^1,\,\ldots,\tx^n),\\
. \ . \ . \ .\ . \ . \ . \ . \ . \ . \ 
. \ . \ . \ . \ .\\
x^n=x^n(\tx^1,\,\ldots,\tx^n),\\
T^{i_1\ldots\,i_r}_{j_1\ldots\,j_s}[P]=\dsize\msum{n}\Sb
h_1,\,\ldots,\,h_r\\
k_1,\,\ldots,\,k_s\endSb S^{i_1}_{h_1}\ldots\,S^{i_r}_{h_r}
\ T^{k_1}_{j_1}\ldots\,T^{k_s}_{j_s}\ \tilde T^{h_1\ldots
\,h_r}_{k_1\ldots\,k_s}[P].
\endcases
\mytag{3.6}
$$
Here again $r=r_P$, $s=s_P$, and $P$ runs from $1$ to $Q$. The 
formulas \mythetag{3.5} and \mythetag{3.6} generalize \mythetag{2.7} 
and \mythetag{2.8} for the case of composite tensor bundles. They 
relate two sets of variables \mythetag{3.3} and \mythetag{3.4}.\par
     The composite tensor bundle $N=T^{r_1\ldots\,r_Q}_{s_1\ldots\,
s_Q}\!M$ is a manifold (like $M$ itself) and \mythetag{3.3} are the 
local coordinates of its point $q$ in some local chart. Therefore, 
one can consider the tangent space $T_q(N)$. By analogy to 
\mythetag{1.5} it is defined as the span of all partial derivatives
with respect to the variables \mythetag{3.3}. Let's denote them
$$
\xalignat 2
&\hskip -2em
\bold U_i=\frac{\partial}{\partial x^i},
&&\bold V^{j_1\ldots\,j_s}_{i_1\ldots\,i_r}[P]=\frac{\partial}
{\partial T^{i_1\ldots\,i_r}_{j_1\ldots\,j_s}[P]
\vphantom{\vrule height 10pt depth 0pt}},
\mytag{3.7}
\endxalignat
$$
where $r=r_P$, $s=s_P$ and $P$ runs from $1$ to $Q$. Note that 
$\bold U_i$ in \mythetag{3.7} are different from $\bold E_i$ in
\mythetag{1.3}. Passing from \mythetag{3.3} to another
set of local coordinates \mythetag{3.4}, one should define the
other set of partial derivatives spanning the tangent space $T_q(N)$:
$$
\xalignat 2
&\hskip -2em
\tilde\bold U_i=\frac{\partial}{\partial\tx^i},
&&\tilde\bold V^{j_1\ldots\,j_s}_{i_1\ldots\,i_r}[P]=\frac{\partial}
{\partial\tilde T^{i_1\ldots\,i_r}_{j_1\ldots\,j_s}[P]
\vphantom{\vrule height 10pt depth 0pt}}.
\mytag{3.8}
\endxalignat
$$
From \mythetag{3.4} and \mythetag{3.5} one easily derives the transformations formulas
$$
\cases
\tilde\bold V^{j_1\ldots\,j_s}_{i_1\ldots\,i_r}[P]=\dsize\msum{n}\Sb
h_1,\,\ldots,\,h_r\\
k_1,\,\ldots,\,k_s\endSb S^{h_1}_{i_1}\ldots\,S^{h_r}_{i_r}
\ T^{j_1}_{k_1}\ldots\,T^{j_s}_{k_s}\ \bold V^{k_1
\ldots\,k_s}_{h_1\ldots\,h_r}[P],\\
\vspace{2ex}
\gathered
\tilde\bold U_j=\sum^n_{i=1}S^i_j\,\bold U_i+\sum^Q_{P=1}
\dsize\msum{n}\Sb i_1,\,\ldots,\,i_r\\j_1,\,\ldots,\,j_s\endSb
\dsize\msum{n}\Sb h_1,\,\ldots,\,h_r\\k_1,\,\ldots,\,k_s\endSb
S^{i_1}_{h_1}\ldots\,S^{i_r}_{h_r}\ T^{k_1}_{j_1}\ldots\,T^{k_s}_{j_s}
\,\times\\
\times\!\left(\,\sum^r_{m=1}\sum^n_{v_m=1}\!\left(\,\shave{\sum^n_{h=1}}
T^{h_m}_h\,\frac{\partial S^h_{v_m}}{\partial\tx^j}\!\right)
\tilde T^{h_1\ldots\,v_m\ldots\,h_r}_{k_1\ldots\,\ldots\,\ldots\,k_s}[P]
\,+\right.\\
\left.+\sum^s_{m=1}\sum^n_{w_m=1}\!
\left(\,\shave{\sum^n_{h=1}}\frac{\partial T^{w_m}_h}
{\partial\tx^j}\,S^h_{k_m}\!\right)\tilde T^{\,h_1\ldots\,\ldots\,\ldots
\,h_r}_{k_1\ldots\,w_m\ldots\,k_s}[P]\!\right)
\bold V^{j_1\ldots\,j_s}_{i_1\ldots\,i_r}[P].
\endgathered
\endcases\
\mytag{3.9}
$$
The transformation formulas \mythetag{3.9} express the tangent vectors
\mythetag{3.8} through the tangent vectors \mythetag{3.7} in $T_q(N)$.
In order to simplify these formulas we introduce the following
$\theta$-parameters defined through the transition matrices
\mythetag{1.2}:
$$
\hskip -2em
\tilde\theta^k_{ij}=\sum^n_{h=1}T^k_h\,\frac{\partial S^h_j}
{\partial\tx^i}=\sum^n_{h=1}T^k_h\,\frac{\partial x^h}
{\partial\tx^i\,\partial\tx^j}.
\mytag{3.10}
$$
Looking at the right hand side of \mythetag{3.10}, we see that
$\tilde\theta^k_{ij}$ are symmetric:
$$
\tilde\theta^k_{ij}=\tilde\theta^k_{\!j\,i}.
$$
Moreover, note that $S$ and $T$ are inverse to each other. From this
fact we derive
$$
\hskip -2em
0=\frac{\partial\delta^k_j}{\partial\tx^i}=\frac{\partial}
{\partial\tx^i}\!\left(\,\shave{\sum^n_{h=1}}T^k_h\,S^h_j
\!\right)=\sum^n_{h=1}T^k_h\frac{\partial S^h_j}{\partial\tx^j}
+\sum^n_{h=1}\frac{\partial T^k_h}{\partial\tx^i}
\,S^h_j.
\mytag{3.11}
$$
Now, comparing \mythetag{3.10} and \mythetag{3.11}, we find that
$\tilde\theta^k_{ij}$ can also be defined as
$$
\hskip -2em
\tilde\theta^k_{ij}=-\sum^n_{h=1}\frac{\partial T^k_h}
{\partial\tx^i}\,S^h_j.
\mytag{3.12}
$$
Applying \mythetag{3.10} and \mythetag{3.12} to the formula \mythetag{3.9},
we can simplify it as follows:
$$
\cases
\tilde\bold V^{j_1\ldots\,j_s}_{i_1\ldots\,i_r}[P]=\dsize\msum{n}\Sb
h_1,\,\ldots,\,h_r\\
k_1,\,\ldots,\,k_s\endSb S^{h_1}_{i_1}\ldots\,S^{h_r}_{i_r}
\ T^{j_1}_{k_1}\ldots\,T^{j_s}_{k_s}\ \bold V^{k_1
\ldots\,k_s}_{h_1\ldots\,h_r}[P],\\
\vspace{2ex}
\gathered
\tilde\bold U_j=\sum^n_{i=1}S^i_j\,\bold U_i+\sum^Q_{P=1}
\dsize\msum{n}\Sb i_1,\,\ldots,\,i_r\\j_1,\,\ldots,\,j_s\endSb
\dsize\msum{n}\Sb h_1,\,\ldots,\,h_r\\k_1,\,\ldots,\,k_s\endSb
\!\left(\,\sum^r_{m=1}\sum^n_{v_m=1}\tilde\theta^{h_m}_{\!j\,v_m}
\times\right.\\
\left.\times\ \tilde T^{h_1\ldots\,v_m\ldots\,h_r}_{k_1\ldots
\,\ldots\,\ldots\,k_s}[P]
-\sum^s_{m=1}\sum^n_{w_m=1}\tilde\theta^{w_m}_{\!j\,k_m}\,
\tilde T^{\,h_1\ldots\,\ldots\,\ldots\,h_r}_{k_1\ldots\,w_m\ldots\,k_s}[P]
\!\right)\!\times\\
\times\ S^{i_1}_{h_1}\ldots\,S^{i_r}_{h_r}\ T^{k_1}_{j_1}
\ldots\,T^{k_s}_{j_s}\ \bold V^{j_1\ldots\,j_s}_{i_1\ldots\,i_r}[P].
\endgathered
\endcases\
\mytag{3.13}
$$
The inverse transformation formulas expressing the tangent vectors
\mythetag{3.7} trough \mythetag{3.8} are written by analogy to
\mythetag{3.13}. They look like
$$
\cases
\bold V^{j_1\ldots\,j_s}_{i_1\ldots\,i_r}[P]=\dsize\msum{n}\Sb
h_1,\,\ldots,\,h_r\\
k_1,\,\ldots,\,k_s\endSb T^{h_1}_{i_1}\ldots\,T^{h_r}_{i_r}
\ S^{j_1}_{k_1}\ldots\,S^{j_s}_{k_s}\ \tilde\bold V^{k_1
\ldots\,k_s}_{h_1\ldots\,h_r}[P],\\
\vspace{2ex}
\gathered
\bold U_j=\sum^n_{i=1}T^i_j\,\tilde\bold U_i+\sum^Q_{P=1}
\dsize\msum{n}\Sb i_1,\,\ldots,\,i_r\\j_1,\,\ldots,\,j_s\endSb
\dsize\msum{n}\Sb h_1,\,\ldots,\,h_r\\k_1,\,\ldots,\,k_s\endSb
\!\left(\,\sum^r_{m=1}\sum^n_{v_m=1}\theta^{h_m}_{\!j\,v_m}
\times\right.\\
\left.\times\ T^{h_1\ldots\,v_m\ldots\,h_r}_{k_1\ldots
\,\ldots\,\ldots\,k_s}[P]
-\sum^s_{m=1}\sum^n_{w_m=1}\theta^{w_m}_{\!j\,k_m}\,
T^{\,h_1\ldots\,\ldots\,\ldots\,h_r}_{k_1\ldots\,w_m\ldots\,k_s}[P]
\!\right)\!\times\\
\times\ T^{i_1}_{h_1}\ldots\,T^{i_r}_{h_r}\ S^{k_1}_{j_1}
\ldots\,S^{k_s}_{j_s}\ \tilde\bold V^{j_1\ldots\,j_s}_{i_1
\ldots\,i_r}[P].
\endgathered
\endcases\
\mytag{3.14}
$$
The $\theta$-parameters in \mythetag{3.14} are defined by analogy to
\mythetag{3.10} and \mythetag{3.12}:
$$
\hskip -2em
\theta^k_{ij}=\sum^n_{h=1}S^k_h\,\frac{\partial T^h_j}
{\partial x^i}=\sum^n_{h=1}S^k_h\,\frac{\partial\tx^h}
{\partial x^i\,\partial x^j}=-\sum^n_{h=1}\frac{\partial S^k_h}
{\partial x^i}\,T^h_j.
\mytag{3.15}
$$
Moreover, we easily derive the following two formulas:
$$
\pagebreak
\align
&\hskip -2em
\theta^k_{ij}=-\sum^n_{h=1}\sum^n_{p=1}\sum^n_{q=1}
\tilde\theta^h_{\!p\,q}S^k_h\,T^p_i\,T^q_j,
\mytag{3.16}\\
&\hskip -2em
\tilde\theta^k_{ij}=-\sum^n_{h=1}\sum^n_{p=1}\sum^n_{q=1}
\theta^h_{\!p\,q}T^k_h\,S^p_i\,S^q_j.
\mytag{3.17}
\endalign
$$
The formulas \mythetag{3.16} and \mythetag{3.17} relate the $\theta$-parameters given by the formula \mythetag{3.14} with
those given by the formulas \mythetag{3.10} and \mythetag{3.12}.
\head
4. Extended tensor fields.
\endhead
\mydefinition{4.1} Let $N=T^{r_1\ldots\,r_Q}_{s_1\ldots\,s_Q}\!M$ 
be a composite tensor bundle over a smooth real manifold $M$. An
extended tensor field $\bold X$ of the type $(\alpha,\beta)$ is a 
tensor-valued function in $N=T^{r_1\ldots\,r_Q}_{s_1\ldots\,s_Q}\!M$
such that it takes each point $q\in N$ to some tensor $\bold X(q)\in T^\alpha_\beta(p,M)$, where $p=\pi(q)$ is the projection of $q$.
\enddefinition
     Informally speaking, an extended tensor field $\bold X$ is 
a tensorial function with one point argument $p\in M$ and $Q$ 
tensorial arguments $\bold T[1],\,\ldots,\,\bold T[Q]$. In a local
chart $U\subset M$ it is represented by its components $X^{i_1\ldots\,
i_\alpha}_{j_1\ldots\,j_\beta}$ each of which is a function of the
variables \mythetag{3.3}. When passing from $U$ to an overlapping
chart $\tilde U$ the components of an extended tensor field are
transformed according to the standard formula \mythetag{1.13}, while
their arguments are transformed according to the formula \mythetag{3.5}.
Tensor bundles, which we considered above, are used to realize 
these arguments geometrically in a coordinate-free form.\par
     The term {\tencyr\char '074}extended tensor field{\tencyr\char '076}
was first introduced in \mycite{11} in order to describe the force field
of a Newtonian dynamical system. Indeed, writing the Newton's second
law for a point mass $m\,\bold a=\bold F(\bold r,\bold v)$, where $\bold v
=\dot\bold r$ is its velocity and $\bold a=\ddot\bold r$ is its
acceleration, we encounter a vector-function $\bold F$ with two arguments
$\bold r$ and $\bold v$. The first argument represents a point of the
$3$-dimensional space $\Bbb E$, while the second argument is a vector
attached to that point. So, both they form a point of the tangent bundle
$T\Bbb E$. The standard $3$-dimensional Euclidean space $\Bbb E$ is a
very simple thing, its tangent bundle $T\Bbb E$ can be treated as a
$6$-dimensional space parametrized by pairs of $3$-dimensional vectors.
However, even in this trivial case, choosing curvilinear coordinates in
$\Bbb E$, we find that the vectors $\bold r$ and $\bold v$ are slightly
different in their nature.\par
     Another place, where we find extended tensor fields with the 
arguments in a tangent bundle $TM$, is the geometry of Finsler (see
\mycite{32}). Here the metric tensor $\bold g$ depends on the velocity
vector $\bold v$. Some generalizations of the Finslerian geometry 
motivated by the Lagrangian dynamics were suggested in \mycite{27} 
and \mycite{28}. A different approach to understanding extended 
tensor fields with the arguments in a tangent bundle $TM$ was used in
\mycite{33}.\par
    Extended tensor fields with the arguments in a cotangent bundle
$T^*\!M$ are natural in Hamiltonian mechanics. Indeed, if a Hamiltonian
dynamical system is produced from a Lagrangian dynamical system, its
Hamilton function $H$ depends on a point of its configuration space
$M$ and on a momentum covector $\bold p$ at this point.\par
    More complicated extended tensor fields can be found in physics of continuous media and in field theories. For example, in \mycite{34} we
find that the specific free energy function $f(T,\hat\bold G)$ depends
on the temperature $T$ and the elastic part of the deformation tensor.
Such a function should certainly be treated as an extended scalar field
with the arguments in $N=T^{00}_{02}\Bbb E$ since the temperature $T$ 
is a scalar field and $\hat\bold G$ is a tensor field of the type $(0,2)$.
\head
5. The algebra of extended tensor fields.
\endhead
    Suppose that some composite tensor bundle $T^{r_1\ldots\,r_Q}_{s_1
\ldots\,s_Q}\!M$ is fixed. Let's denote by $T^\alpha_\beta(M)$ the set
of all extended tensor fields of the type $(\alpha,\beta)$. According
to the above definition~\mythedefinition{4.1} an extended tensor field
is a function, its values are tensors. Hence, we can perform the
algebraic operations with these values and thus define the algebraic
operations with extended tensor fields:
\roster
\item $T^\alpha_\beta(M)+T^\alpha_\beta(M)\longrightarrow
T^\alpha_\beta(M)$;
\item $T^0_0(M)\otimes T^\alpha_\beta(M)\longrightarrow
T^\alpha_\beta(M)$;
\item $T^\alpha_\beta(M)\otimes T^\sigma_\mu(M)\longrightarrow
T^{\alpha+\sigma}_{\beta+\mu}(M)$;
\endroster
According to the definition~\mythedefinition{4.1}, each extended
scalar field is a real-valued function in the bundle $N$. Such
functions form a ring, we denote it $\goth F(N)$. Then
$$
\hskip -2em
T^0_0(M)=\goth F(N).
\mytag{5.1}
$$
Due to the properties \therosteritem{1} and \therosteritem{2} the 
set of extended tensor fields $T^\alpha_\beta(M)$ of the fixed type
$(\alpha,\beta)$ is a module over the ring \mythetag{5.1}.\par
     In general, one cannot add tensor fields of two different types $(\alpha,\beta)\neq (\sigma,\mu)$. However, one can consider formal 
sums
$$
\hskip -2em
\bold X=\bold X[1]+\bold X[2]+\ldots+\bold X[K],
\mytag{5.2}
$$
where $\bold X[1],\ \ldots,\ \bold X[K]$ are taken from various modules 
$T^\alpha_\beta(M)$ over the ring $\goth F(N)$. By definition, the set of
all sums \mythetag{5.2} is called the {\it direct sum}
$$
\hskip -2em
\bold T(M)=\bigoplus^\infty_{\alpha=0}
\bigoplus^\infty_{\beta=0}T^\alpha_\beta(M).
\mytag{5.3}
$$
With respect to the algebraic operations \therosteritem{1},
\therosteritem{2}, and \therosteritem{3} the direct sum \mythetag{5.3}
is a {\it graded algebra} over the ring $T^0_0(M)$. This algebra is called
the {\it algebra of extended tensor fields\/} (or the {\it extended 
algebra} for short).\par
     The fourth class of algebraic operations in the extended algebra
\mythetag{5.3} is formed by the operations of contraction. They are
performed with respect to some pair of indices one of which is an upper
index and the other is a lower index (see details in \mycite{1} and
\mycite{2}). For the sake of simplicity we denote these operations
as follows:
\roster
\item[4] $C\!:\,T^\alpha_\beta(M)\longrightarrow T^{\alpha-1}_{\beta-1}(M)$
for $\alpha\geqslant 1$ and $\beta\geqslant 1$.
\endroster
As it follows from the property \therosteritem{4}, the contraction operations are concordant with the structure of graded algebra in 
$\bold T(M)$. The same is true for the operations of addition and 
tensor multiplication.\par
\head
6. Differentiation of tensor fields.
\endhead
\mydefinition{6.1}  An extended tensor field $\bold X$ of the type $(\alpha,
\beta)$ associated with a composite tensor bundle $N=T^{r_1\ldots\,r_Q}_{s_1
\ldots\,s_Q}\!M$ is called {\it smooth\/} if its components $X^{i_1\ldots\,
i_\alpha}_{j_1\ldots\,j_\beta}$ are smooth function of their arguments
\mythetag{3.3} for any local chart $U$ of $M$.
\enddefinition
    From now on we shall restrict our previous definition \mythetag{5.3} 
of the extended algebra and denote by $\bold T(M)$ the algebra of smooth
extended tensor fields. Similarly, by $\goth F(N)$ below we denote the 
ring of smooth scalar functions in $N$. Then the equality \mythetag{5.1}
remains valid. The basic idea of considering smooth smooth extended tensor fields is to introduce the operation of differentiation in addition to the above four algebraic operations in $\bold T(M)$. 
\mydefinition{6.2} A mapping $D\!:\,\bold T(M)\to\bold T(M)$ is called 
a {\it differentiation} of the extended algebra of tensor fields if the
following conditions are fulfilled:
\roster
\rosteritemwd=10pt
\item $D$ is concordant with the grading: $D(T^\alpha_\beta(M))
\subset T^\alpha_\beta(M)$;
\item $D$ is $\Bbb R$-linear: \hskip -0.95cm\vtop{\hsize 6.4cm
      \noindent $D(\bold X+\bold Y)=D(\bold X)+D(\bold Y)$
      \newline and $D(\lambda\bold X)=\lambda D(\bold X)$
      for $\lambda\in\Bbb R$;\vskip 1.3ex}
\item $D$ commutates with the contractions: $D(C(\bold X))=C(D(\bold X))$;
\item $D$ obeys the Leibniz rule: $D(\bold X\otimes\bold Y)=D(\bold X)
      \otimes\bold Y+\bold X\otimes D(\bold Y)$.
\endroster
\enddefinition
\noindent Let's consider the set of all differentiations of the extended
algebra $\bold T(M)$. We denote it $\goth D(M)$. It is easy to check up
that 
\roster
\rosteritemwd=10pt
\item the sum of two differentiations is a differentiation of the algebra
$\bold T(M)$;
\item the product of a differentiation by a smooth function in $N$ is a differentiation of the algebra $\bold T(M)$.
\endroster
Now we see that $\goth D(M)$ is equipped with the structure 
of a module over the ring of smooth functions $\goth F(N)=T^0_0(M)$. 
The composition of  two differentiations $D_1$ and $D_2$ is not a
differentiation, but their commutator 
$$
\hskip -2em
[D_1,\,D_2]=D_1\compos D_2-D_2\compos D_1
\mytag{6.1}
$$
is a differentiation. Therefore, $\goth D(M)$ is a Lie algebra. Note,
however, that $\goth D(M)$ is not a Lie algebra over the ring of smooth
functions $\goth F(N)$. It is only a Lie algebra over the field of real
numbers $\Bbb R$. For this reason the structure of Lie algebra in
$\goth D(M)$ is not of a primordial importance.\par
\head
7. Localization.
\endhead
    Smooth extended tensor fields are global objects related to the
tensor bundle $N$ in whole, but they are functions --- their values
are local objects so that two different tensor fields $\bold A\neq 
\bold B$ can take the same values at some particular points. Whenever
this happens, we write $\bold A_q=\bold B_q$, where $q\in N$.\par
    Differentiations of the algebra $\bold T(M)$, as they are introduced
in the definition~\mythedefinition{6.2}, are global objects without any
explicit subdivision into parts related to separate points of the bundle
$N$. Below in this section we shall show that they also can be represented
as functions taking their values in some linear spaces associated with
separate points of the manifold $N$.
    Let $D\in\goth D(M)$ be a differentiation of the algebra of extended
tensor fields $\bold T(M)$. Let's denote by $\delta$ the restriction of 
the mapping $D\!:\bold T(M)\to\bold T(M)$ to the module $T^0_0(M)$ 
of extended scalar fields in \mythetag{5.3}:
$$
\hskip -2em
\delta\!:\,T^0_0(M)\to T^0_0(M).
\mytag{7.1}
$$
Since $T^0_0(M)=\goth F(M)$, the mapping $\delta$ in \mythetag{7.1} 
is a differentiation of the ring of smooth functions in $N$. It is 
known (see \S\,1 in Chapter~\uppercase\expandafter{\romannumeral 1} 
of \mycite{35}) that any differentiation of the ring of smooth functions 
of an arbitrary smooth manifold is determined by some vector field 
$\bold Z$ in this manifold. Applying this fact to $N$, we get
$$
\pagebreak
\hskip -2em
\delta=\bold Z=\sum^n_{i=1}Z^i\ \bold U_i+\sum^Q_{P=1}
\dsize\msum{n}\Sb i_1,\,\ldots,\,i_r\\j_1,\,\ldots,\,j_s\endSb
Z^{i_1\ldots\,i_r}_{j_1\ldots\,j_s}[P]\ \bold V^{j_1\ldots\,j_s}_{i_1
\ldots\,i_r}[P],
\mytag{7.2}
$$
where $\bold U_i$ and $\bold V^{j_1\ldots\,j_s}_{i_1\ldots\,i_r}[P]$ 
are the differential operators \mythetag{3.7}. From the representation
\mythetag{7.2} for the mapping \mythetag{7.1} we immediately derive
the following lemma.
\mylemma{7.1} Let $\psi\in T^0_0(M)=\goth F(N)$ an extended scalar 
field (a smooth function) identically constant within some open
subset $O\subset N$ and let $\varphi=D(\psi)$ for some differentiation
$D\in\goth D(M)$. Then $\varphi=0$ within the open set $O$.
\endproclaim
\demo{Proof} Since $D(\psi)=\delta(\psi)$, choosing some local chart
$U$ and applying the differential operators \mythetag{3.7} to a constant,
from \mythetag{7.2} we derive that $\delta(\psi)=\bold Z\psi=0$ at any
point $q$ of the open set $O$.\qed\enddemo
\mylemma{7.2} Let\/ $\bold X$ be an extended tensor field of the type
$(\alpha,\beta)$. If\/ $\bold X\equiv 0$, then $D(\bold X)$ is also
identically equal to zero for any differentiation $D\in\goth D(M)$.
\endproclaim
    The proof is trivial. Since $\bold X\equiv 0$, we can write $\bold
X=\lambda\,\bold X$ with $\lambda\neq 1$. Then, applying the item
\therosteritem{2} of the definition~\mythedefinition{6.2}, we get 
$D(\bold X)=\lambda\,D(\bold X)$. Since $\lambda\neq 1$, this yields 
the required equality $D(\bold X)\equiv 0$.
\mylemma{7.3} Let $\bold X$ be an extended tensor field of the type
$(\alpha,\beta)$ identically zero within some open set $O\subset N$.
If\/ $\bold Y=D(\bold X)$ for some differentiation $D\in\goth D(M)$,
then $\bold Y_{\!q}=0$ at any point $q\in O$.
\endproclaim
\demo{Proof} Let's choose some arbitrary point $q\in O$ and take some
smooth scalar function $\eta\in\goth F(N)$ such that it is identically
equal to the unity in some open neighborhood $O{\,'}\subset O$ of the point
$q$ and identically equal to zero outside the open set $O$. The existence
of such a function $\eta$ is easily proved by choosing some local chart $U$
that covers the point $q$. The product $\eta\,\bold X$ is identically equal
to zero:
$$
\hskip -2em
\eta\otimes\bold X=\eta\,\bold X\equiv 0.
\mytag{7.3}
$$
Applying the differentiation $D$ to \mythetag{7.3}, then taking into 
account the lemma~\mythelemma{7.2} and the item \therosteritem{4} of 
the definition~\mythedefinition{6.2}, we obtain
$$
0=D(0)=D(\eta\otimes\bold X)=D(\eta)\otimes\bold X+\eta\otimes D(\bold X)
=D(\eta)\,\bold X+\eta\,D(\bold X).
$$
Note that $D(\eta)=0$ at the point $q$ due to the lemma~\mythelemma{7.1}.
Moreover, $\bold X_q=0$ and $\eta=1$ at the point $q$. Therefore, by
specifying the above equality to the point $q$ we get $D(\bold X)=0$ at
the point $q$. The lemma is proved.\qed\enddemo
\mylemma{7.4} If two extended tensor fields $\bold X$ and $\bold Y$ 
are equal within some open neighborhood $O$ of a point $q\in N$, then 
for any differentiation $D\in\goth D(M)$ their images $D(\bold X)$ and
$D(\bold Y)$ are equal at the point $q$.
\endproclaim
    The lemma~\mythelemma{7.4} follows immediately from the 
lemma~\mythelemma{7.3}. This lemma is a basic tool for our purposes
of localization in this section.\par
    Let $q$ be some point of $N$ and let $\pi(q)$ be its projection
in the base manifold $M$. Taking some local chart $U$ that covers the
point $p$ in $M$, we can use its preimage $\pi^{-1}(U)$ as a local
chart in $N$ covering the point $q$. The variables \mythetag{3.3} 
form the complete set of local coordinates in the chart $\pi^{-1}(U)$.
Any extended tensor field $\bold X$ of the type $(\alpha,\beta)$ is
represented  by the formula 
$$
\hskip -2em
\bold X=\sum^n_{i_1=1}\ldots\sum^n_{i_\alpha=1}\sum^n_{j_1=1}\ldots
\sum^n_{j_\beta=1}X^{i_1\ldots\,i_\alpha}_{j_1\ldots\,j_\beta}\ 
\bold E^{j_1\ldots\,j_\beta}_{i_1\ldots\,i_\alpha},
\mytag{7.4}
$$
which is identical to \mythetag{1.11}. The only difference is that 
the coefficients $X^{i_1\ldots\,i_\alpha}_{j_1\ldots\,j_\beta}$ in
\mythetag{7.4} depend not only on $x^1,\,\ldots,x^n$, but on the 
whole set of variables \mythetag{3.3}. Taking some differentiation
$D\in\goth D(M)$, we can apply it to the left hand side of the equality
\mythetag{7.4}, but we cannot apply it to each summand in the right hand
side of \mythetag{7.4}. The matter is that the scalars $X^{i_1\ldots\,
i_\alpha}_{j_1\ldots\,j_\beta}$ and the tensors $\bold E^{j_1\ldots\,
j_\beta}_{i_1\ldots\,i_\alpha}$ are defined locally only within the chart
$\pi^{-1}(U)$, therefore they do not fit the
definition~\mythedefinition{4.1}. In order to convert them to global 
fields we choose some smooth function $\eta\in\goth F(N)$ such that it is
identically equal to the unity within some open neighborhood of the point
$q$ and is identically zero outside the chart $\pi^{-1}(U)$. Let's define
$$
\align
\hskip -20em
\hat X^{i_1\ldots\,i_\alpha}_{j_1\ldots\,j_\beta}&=\cases
\eta\ X^{i_1\ldots\,i_\alpha}_{j_1\ldots\,j_\beta}&\text{within \ }
\pi^{-1}(U),\\
\ \ 0&\text{outside \ }\pi^{-1}(U),\endcases
\mytag{7.5}\\
\vspace{2ex}
\hskip -2em
\hat\bold E_i&=\cases\eta\ \bold E_i &\text{within \ }
\pi^{-1}(U),\\
\ \ 0&\text{outside \ }\pi^{-1}(U),\endcases
\mytag{7.6}\\
\vspace{2ex}
\hskip -2em
\hbh^i&=\cases\eta\ dx^i &\text{within \ }
\pi^{-1}(U),\\
\ \ 0&\text{outside \ }\pi^{-1}(U).\endcases
\mytag{7.7}
\endalign
$$
Then by analogy to \mythetag{1.10} we write
$$
\hskip -2em
\hat\bold E^{j_1\ldots\,j_\beta}_{i_1\ldots\,i_\alpha}
=\hat\bold E_{i_1}\otimes\ldots\otimes\hat\bold E_{i_\alpha}
\otimes \hbh^{j_1}\otimes\ldots\otimes\hbh^{j_\beta}.
\mytag{7.8}
$$
Taking into account \mythetag{7.5}, \mythetag{7.6}, \mythetag{7.7},
and \mythetag{7.8}, from \mythetag{7.4} we derive
$$
\hskip -2em
\eta^m\ \bold X=\sum^n_{i_1=1}\ldots\sum^n_{i_\alpha=1}\sum^n_{j_1=1}
\ldots\sum^n_{j_\beta=1}\hat X^{i_1\ldots\,i_\alpha}_{j_1\ldots\,j_\beta}\ 
\hat\bold E^{j_1\ldots\,j_\beta}_{i_1\ldots\,i_\alpha},
\mytag{7.9}
$$
where $m=\alpha+\beta+1$. Now we can apply $D$ to each summand in the right
hand side of the equality \mythetag{7.9}. Using the item \therosteritem{4}
of the definition~\mythedefinition{6.2}, we get
$$
\hskip -2em
\gathered
D(\eta^m\ \bold X)=\sum^n_{i_1=1}\ldots\sum^n_{i_\alpha=1}
\sum^n_{j_1=1}\ldots\sum^n_{j_\beta=1}D(\hat X^{i_1\ldots\,
i_\alpha}_{j_1\ldots\,j_\beta})\ \hat\bold E^{j_1\ldots\,
j_\beta}_{i_1\ldots\,i_\alpha}\,+\\
+\,\sum^n_{i_1=1}\ldots\sum^n_{i_\alpha=1}\sum^n_{j_1=1}
\ldots\sum^n_{j_\beta=1}\hat X^{i_1\ldots\,i_\alpha}_{j_1
\ldots\,j_\beta}\ D(\hat\bold E^{j_1\ldots\,j_\beta}_{i_1
\ldots\,i_\alpha}).
\endgathered
\mytag{7.10}
$$
Due to the lemma~\mythelemma{7.4} we have $D(\eta^m\ \bold X)
=D(\bold X)$ at the point $q$. Moreover, due to \mythetag{7.5},
\mythetag{7.6}, \mythetag{7.7}, \mythetag{7.8} and since $\eta(q)=1$,
we have $\hat\bold E^{j_1\ldots\,j_\beta}_{i_1\ldots\,i_\alpha}
=\bold E^{j_1\ldots\,j_\beta}_{i_1\ldots\,i_\alpha}$ and $\hat X^{i_1\ldots\,i_\alpha}_{j_1\ldots\,j_\beta}=X^{i_1\ldots\,i_\alpha}_{j_1\ldots\,j_\beta}$ at this point. As for the fields $D(\hat X^{i_1\ldots\,i_\alpha}_{j_1\ldots\,j_\beta})$ and $D(\hat\bold E^{j_1
\ldots\,j_\beta}_{i_1\ldots\,i_\alpha})$ in \mythetag{7.10}, again 
due to the lemma~\mythelemma{7.4} their values at the point $q$ do not
depend on a particular choice of the function $\eta$. Since $D(\hat
X^{i_1\ldots\,i_\alpha}_{j_1\ldots\,j_\beta})=\delta(\hat X^{i_1\ldots\,
i_\alpha}_{j_1\ldots\,j_\beta})$, from the formulas \mythetag{7.2} and
\mythetag{3.7} for the value of $D(\hat X^{i_1\ldots\,i_\alpha}_{j_1
\ldots\,j_\beta})$ at the point $q$ we derive
$$
\pagebreak
D(\hat X^{i_1\ldots\,i_r}_{j_1\ldots\,j_s})=\sum^n_{i=1}Z^i\ 
\frac{\partial X^{i_1\ldots\,i_\alpha}_{j_1\ldots\,j_\beta}}
{\partial x^i}+\sum^Q_{P=1}Z^{i_1\ldots\,i_r}_{j_1\ldots\,j_s}[P]
\ \frac{\partial X^{i_1\ldots\,i_\alpha}_{j_1\ldots\,j_\beta}}
{\partial T^{i_1\ldots\,i_r}_{j_1\ldots\,j_s}[P]
\vphantom{\vrule height 10pt depth 0pt}}.\quad
\mytag{7.11}
$$
In order to evaluate $D(\hat\bold E^{j_1\ldots\,j_\beta}_{i_1\ldots\,
i_\alpha})$ at the point $q$ let's apply $D$ to \mythetag{7.8}. Using 
the item \therosteritem{4} of the definition~\mythedefinition{6.2}, 
we obtain the following equality:
$$
\hskip -2em
\gathered
D(\hat\bold E^{j_1\ldots\,j_\beta}_{i_1\ldots\,i_\alpha})=\sum^\alpha_{v=1}
\bold E_{i_1}\otimes\ldots\otimes D(\,\hat\bold E_{i_v})\otimes
\ldots\otimes\bold E_{i_\alpha}\otimes\\
\otimes\ dx^{j_1}\otimes\ldots\otimes dx^{j_\beta}
+\sum^\beta_{w=1}\bold E_{i_1}\otimes\ldots\otimes\bold E_{i_\alpha}
\otimes\\
\otimes\ dx^{j_1}\otimes\ldots\otimes D(\,\hbh^{j_w})\otimes\ldots
\otimes dx^{j_\beta}.\vphantom{\sum^s_{w=1}}
\endgathered
\mytag{7.12}
$$
Summarizing the above three formulas \mythetag{7.10}, \mythetag{7.11}, 
\mythetag{7.12}, and the equality $D(\eta^m\,\bold X)=D(\bold X)$,
we can formulate the following lemma.
\mylemma{7.5} Any differentiation $D\in\goth D(M)$ is uniquely fixed by 
its restrictions to the modules $T^0_0(M)$, $T^1_0(M)$, $T^0_1(M)$ in
the direct sum \mythetag{5.3}.
\endproclaim
     The value of the extended vector field $D(\hat\bold E_i)$ at the
the point $q$ is a vector of the tangent space $T_{\pi(q)}(M)$. We can
write the following expansion for this vector:
$$
\hskip -2em
D(\hat\bold E_i)=\sum^n_{k=1}\Gamma^k_i\,\bold E_k.
\mytag{7.13}
$$
Due to the lemma~\mythelemma{7.4} the left hand side of the equality
\mythetag{7.13} does not depend on a particular choice of the function
$\eta$ in \mythetag{7.6}. Therefore, the coefficients $\Gamma^k_i$ in
\mythetag{7.13} represent the differentiation $D$ at the point $q$ for
a given local chart $U$ in $M$. The same is true for $Z^i$ and
$Z^{i_1\ldots\,i_r}_{j_1\ldots\,j_s}[P]$ in \mythetag{7.11}. Being 
dependent on $q$, all these coefficients $Z^i$, $Z^{i_1\ldots\,i_r}_{j_1\ldots\,j_s}[P]$, and $\Gamma^k_i$ are some 
smooth functions of the variables \mythetag{3.3}. However, if $q$ is 
fixed, they all are constants.\par
    Under a change of a local chart the coefficients $Z^i$,
$Z^{i_1\ldots\,i_r}_{j_1\ldots\,j_s}[P]$, and $\Gamma^k_i$ obey some
definite transformation rules. The transformation rules for $Z^i$
and $Z^{i_1\ldots\,i_r}_{j_1\ldots\,j_s}[P]$ are derived from the
formulas \mythetag{3.13}, \mythetag{3.14}, and \mythetag{7.2}:
$$
\cases
Z^i=\dsize\sum^n_{j=1}S^i_j\,\tilde Z^j,\\
\vspace{2ex}
\gathered
Z^{i_1\ldots\,i_r}_{j_1\ldots\,j_s}[P]=\dsize\msum{n}\Sb
h_1,\,\ldots,\,h_r\\
k_1,\,\ldots,\,k_s\endSb 
S^{i_1}_{h_1}\ldots\,S^{i_r}_{h_r}
\ T^{k_1}_{j_1}\ldots\,T^{k_s}_{j_s}\ 
\tilde Z^{h_1\ldots\,h_r}_{k_1\ldots\,k_s}[P]\ -\\
-\sum^r_{m=1}\sum^n_{i=1}\sum^n_{j=1}\sum^n_{v_m=1}
\theta^{i_m}_{i\,v_m}\,T^{\,i_1\ldots\,v_m\ldots\,i_r}_{j_1
\ldots\,\ldots\,\ldots\,j_s}[P]\ S^i_j\ \tilde Z^j\ +\\
+\sum^s_{m=1}\sum^n_{i=1}\sum^n_{j=1}\sum^n_{w_m=1}
\theta^{w_m}_{i\,j_m}\,T^{\,i_1\ldots\,\ldots\,\ldots\,
i_r}_{j_1\ldots\,w_m\ldots\,j_s}[P]\ S^i_j\ \tilde Z^j.
\endgathered
\endcases\qquad
\mytag{7.14}
$$
Here the components of the direct and inverse transition matrices $S$
and $T$ are taken from \mythetag{1.2}, while the $\theta$-parameters 
are given by the formulas \mythetag{3.15}. The transformation formulas
\mythetag{7.14} can be inverted in the following way:
$$
\cases
\tilde Z^i=\dsize\sum^n_{j=1}T^i_j\,Z^j,\\
\vspace{2ex}
\gathered
\tilde Z^{\,i_1\ldots\,i_r}_{j_1\ldots\,j_s}[P]=\dsize\msum{n}\Sb
h_1,\,\ldots,\,h_r\\
k_1,\,\ldots,\,k_s\endSb 
T^{i_1}_{h_1}\ldots\,T^{i_r}_{h_r}
\ S^{k_1}_{j_1}\ldots\,S^{k_s}_{j_s}\ 
Z^{h_1\ldots\,h_r}_{k_1\ldots\,k_s}[P]\ -\\
-\sum^r_{m=1}\sum^n_{i=1}\sum^n_{j=1}\sum^n_{v_m=1}
\tilde\theta^{i_m}_{i\,v_m}\,\tilde T^{\,i_1\ldots\,v_m\ldots\,i_r}_{j_1
\ldots\,\ldots\,\ldots\,j_s}[P]\ T^i_j\ Z^j\ +\\
+\sum^s_{m=1}\sum^n_{i=1}\sum^n_{j=1}\sum^n_{w_m=1}
\tilde\theta^{w_m}_{i\,j_m}\,\tilde T^{\,i_1\ldots\,\ldots\,\ldots\,
i_r}_{j_1\ldots\,w_m\ldots\,j_s}[P]\ T^i_j\ Z^j.
\endgathered
\endcases\qquad
\mytag{7.15}
$$
The $\theta$-parameters for \mythetag{7.15} are taken from \mythetag{3.10}
or from \mythetag{3.12}. The transformation formulas \mythetag{7.14} and
\mythetag{7.15} should be completed with the analogous formulas for the
coefficients $\Gamma^k_i$ in the expansion \mythetag{7.13}:
$$
\align
&\hskip -2em
\Gamma^k_i=\sum^n_{b=1}\sum^n_{a=1}S^k_a\,T^b_i\ \tilde\Gamma^a_b+
\sum^n_{a=1}Z^a\,\theta^k_{ai},
\mytag{7.16}\\
&\hskip -2em
\tilde\Gamma^k_i=\sum^n_{b=1}\sum^n_{a=1}T^k_a\,S^b_i\ \Gamma^a_b+
\sum^n_{a=1}\tilde Z^a\,\tilde\theta^k_{ai}.
\mytag{7.17}
\endalign
$$
The formula \mythetag{7.16} is derived immediately from \mythetag{7.13}
and \mythetag{1.4}, taking into account \mythetag{3.15}. The formula
\mythetag{7.17} then is written by analogy.\par
     Let's return back to \mythetag{7.12}. In order to calculate 
$\bold D(\,\hbh_j)$ in this formula we need to remember \mythetag{1.9}.
This equality can be written as follows:
$$
\hskip -2em
C(dx^i\otimes\bold E_j)=\delta^i_j.
\mytag{7.18}
$$
Multiplying the equality \mythetag{7.18} by $\eta^2$, we transform 
it to the following one:
$$
\hskip -2em
C(\hbh^i\otimes\hat\bold E_j)=\delta^i_j\,\eta^2.
\mytag{7.19}
$$
Here $\eta$ is the same function as in \mythetag{7.6} and \mythetag{7.7}.
Now one can apply the differentiation $D$ to both sides of \mythetag{7.19}.
Taking into account the item \therosteritem{3} in the 
definition~\mythedefinition{6.2} and taking into account the
lemma~\mythelemma{7.4}, at the fixed point $q$ we obtain
$$
\hskip -2em
C(D(\hbh^i)\otimes\hat\bold E_j)+C(\hbh^i\otimes D(\hat\bold E_j))
=D(\delta^i_j\,\eta^2)=0.
\mytag{7.20}
$$
Since $\eta=1$ at the point $q$, from \mythetag{7.13}, \mythetag{7.19},
and \mythetag{7.20} we derive
$$
\hskip -2em
C(D(\hbh^i)\otimes\hat\bold E_j)=-\sum^n_{k=1}\Gamma^k_j\ C(\hbh^i
\otimes\hat\bold E_k)=-\Gamma^i_j.
\mytag{7.21}
$$
By the definition of a differentiation the value of $D(\hbh^i)$ \pagebreak 
at the point $q$ is a covector from the cotangent space $T^*_q(M)$. One can
expand it in the basis of differentials $dx^1,\,\ldots,\,dx^n$. Due to
\mythetag{7.21} this expansion looks like
$$
\hskip -2em
D(\hbh^i)=-\sum^n_{j=1}\Gamma^i_j\,dx^j.
\mytag{7.22}
$$
The formula \mythetag{7.22} means that we can strengthen the 
lemma~\mythelemma{7.5} as follows.
\mytheorem{7.1} Any differentiation $D\in\goth D(M)$ is uniquely fixed 
by its restrictions to the modules $T^0_0(M)$, $T^1_0(M)$ in the direct 
sum \mythetag{5.3}.
\endproclaim
    Indeed, let's denote by $\xi$ and $\zeta$ the restrictions of the differentiation $D$ to the modules $T^1_0(M)$ and $T^0_1(M)$ respectively.
Its restriction to the module $T^0_0(M)$ was already considered. It was
denoted by $\delta$ (see formula \mythetag{7.1} above):
$$
\align
&\hskip -2em
\xi\!:\,T^1_0(M)\to T^1_0(M),
\mytag{7.23}\\
&\hskip -2em
\zeta\!:\,T^0_1(M)\to T^0_1(M).
\mytag{7.24}
\endalign
$$
The formulas \mythetag{7.13} and \mythetag{7.22} mean that the mapping
\mythetag{7.24} is completely determined if the mapping \mythetag{7.23}
is known. Hence, $D$ is completely determined by the mappings $\delta$ 
and $\xi$. Then, using \mythetag{7.10}, \mythetag{7.11}, and 
\mythetag{7.12}, one can reconstruct the mapping $D$ itself. Let
$\bold Y=D(\bold X)$ and assume that $\bold X$ is given by the formula
\mythetag{7.4} in local coordinates. Then $\bold Y$ is given by a similar
expansion
$$
\hskip -2em
\bold Y=\sum^n_{i_1=1}\ldots\sum^n_{i_\alpha=1}\sum^n_{j_1=1}\ldots
\sum^n_{j_\beta=1}Y^{i_1\ldots\,i_\alpha}_{j_1\ldots\,j_\beta}\ 
\bold E^{j_1\ldots\,j_\beta}_{i_1\ldots\,i_\alpha},
\mytag{7.25}
$$
where its components $Y^{i_1\ldots\,i_\alpha}_{j_1\ldots\,j_\beta}$ are
calculated as follows:
$$
\gathered
Y^{i_1\ldots\,i_\alpha}_{j_1\ldots\,j_\beta}=
\sum^n_{i=1}Z^i\ \frac{\partial X^{i_1\ldots\,i_\alpha}_{j_1\ldots\,
j_\beta}}{\partial x^i}+\sum^Q_{P=1}\dsize\msum{n}\Sb h_1,\,\ldots,\,h_r\\
k_1,\,\ldots,\,k_s\endSb Z^{h_1\ldots\,h_r}_{k_1\ldots\,k_s}[P]
\ \frac{\partial X^{i_1\ldots\,i_\alpha}_{j_1\ldots\,j_\beta}}
{\partial T^{h_1\ldots\,h_r}_{k_1\ldots\,k_s}[P]
\vphantom{\vrule height 10pt depth 0pt}}\ +\\
+\sum^\alpha_{m=1}\sum^n_{v_m=1}\Gamma^{i_m}_{v_m}\ X^{\,i_1\ldots\,v_m
\ldots\,i_\alpha}_{j_1\ldots\,\ldots\,\ldots\,j_\beta}-\sum^\beta_{m=1}
\sum^n_{w_m=1}\Gamma^{w_m}_{j_m}\ X^{\,i_1\ldots\,\ldots\,\ldots\,
i_\alpha}_{j_1\ldots\,w_m\ldots\,j_\beta}.
\endgathered\qquad
\mytag{7.26}
$$
The formulas \mythetag{7.25} and \mythetag{7.26} prove the 
theorem~\mythetheorem{7.1}. As for the mappings \mythetag{7.23} and
\mythetag{7.24}, one can formulate the following theorem for them.
\mytheorem{7.2} Defining a differentiation $D$ of the algebra of 
extended tensor fields $\bold T(M)$ is equivalent to defining two 
$\Bbb R$-linear mappings \mythetag{7.23} and \mythetag{7.24} such that
$$
\align
&\hskip -6em
\delta(\varphi\,\psi)=\delta(\varphi)\,\psi+\varphi\,\delta(\psi)
\text{ \ for any \ }\varphi,\psi\in\goth F(N)=T^0_0(M),
\mytag{7.27}\\
&\hskip -6em
\xi(\varphi\,\bold X)=\delta(\varphi)\ \bold X+
\varphi\ \xi(\bold X)\text{ \ for any \ }\varphi\in\goth F(N)
\text{ \ and \ }\bold X\in T^1_0(M).\hskip -2em
\mytag{7.28}
\endalign
$$
\endproclaim
The formula \mythetag{7.27} is immediate from the 
definition~\mythedefinition{6.2}. It means that $\delta$ is a
differentiation of the ring $\goth F(N)$. This fact was used for
writing \mythetag{7.2}. The formula \mythetag{7.28} is also immediate
from the definition~\mythedefinition{6.2}.\par
    Looking at the formula \mythetag{7.26}, we see that the differentiation
$D$ acts as a first order linear differential operator upon the components
\pagebreak of a tensor field $\bold X$. The coefficients $Z^i$, $Z^{h_1
\ldots\,h_r}_{k_1\ldots\,k_s}[P]$, $\Gamma^k_i$ of this linear operator are not differentiated in \mythetag{7.26}. Therefore, fixing some point $q\in N$
and taking the values of the coefficients $Z^i$, $Z^{h_1\ldots\,h_r}_{k_1
\ldots\,k_s}[P]$, $\Gamma^k_i$, we can say that we know this linear 
operator at that particular point $q$. In order to formalize this idea we
need a coordinate-free definition of a first order linear differential
operator at a point $q$ acting upon tensors.
\mydefinition{7.1} Two smooth extended tensor fields $\bold X_1$ and 
$\bold X_2$ defined in some open neighborhoods $O_1$ and $O_2$ of a 
point $q\in N$ are called {\it $q$-equivalent\/} if they take the same 
values within some smaller neighborhood $O\subset O_1\cap\O_2$ of the 
point $q$.
\enddefinition
\mydefinition{7.2} A class of $q$-equivalent smooth extended tensor fields
is called a {\it stalk\/} of a smooth extended tensor field at the point $q$.
\enddefinition
The stalks of various extended tensor fields at a fixed point $q\in N$ form
a graded algebra over the real numbers $\Bbb R$ (compare to \mythetag{5.3}
above):
$$
\goth T(q,M)=\bigoplus^\infty_{\alpha=0}\bigoplus^\infty_{\beta=0}
\goth T^\alpha_\beta(q,M).
\mytag{7.29}
$$
Each summand $\goth T^\alpha_\beta(q,M)$ in \mythetag{7.29} is a linear
space over $\Bbb R$. It is composed by stalks of extended tensor fields 
of some fixed type $(\alpha,\beta)$.\par
    A stalk of a tensor field is somewhat like its value at a fixed point
$q$. However, they do not coincide since the stalk comprises much more
information:
$$
\goth T^\alpha_\beta(q,M)\neq T^\alpha_\beta(q,M).
$$
Informally speaking, a stalk is the restriction of a tensor field to the
infinitesimal neighborhood of a point $q$. The information stored in a 
stalk of a tensor field is sufficient for to apply a differential operator
to it. As a result we get a tensor (not a stalk) at a fixed point.
\mydefinition{7.3} A {\it tensorial first order differential operator} 
$D_q$ at a point $q\in N$ is a mapping $D_q\!:\,\goth T(q,M)\to T(q,M)$
such that
\roster
\rosteritemwd=10pt
\item $D_q$ is concordant with the grading: $D_q(\goth T^\alpha_\beta(q,M))
      \subset T^\alpha_\beta(q,M)$;
\item $D_q$ is $\Bbb R$-linear: \hskip -0.95cm\vtop{\hsize 6.4cm
      \noindent $D_q(\bold X+\bold Y)=D_q(\bold X)+D_q(\bold Y)$
      \newline and $D_q(\lambda\,\bold X)=\lambda\,D_q(\bold X)$
      for $\lambda\in\Bbb R$;\vskip 1.3ex}
\item $D_q$ commutates with the contractions: $D_q(C(\bold X))
      =C(D_q(\bold X))$;
\item $D_q$ obeys the Leibniz rule: $D_q(\bold X\otimes\bold Y)
      =D_q(\bold X)\otimes\bold Y_{\!q}+\bold X_q\otimes D_q(\bold Y)$.
\endroster
\enddefinition
The formulas \mythetag{7.25} and \mythetag{7.26} yield the representation
of a tensorial differential operator in a local chart. In the case of a
differential operator $D_q$ the coefficients $Z^i$, $Z^{h_1\ldots\,h_r}_{k_1
\ldots\,k_s}[P]$, $\Gamma^k_i$ are constants related to a point $q$. They
are called the {\it components\/} of $D_q$. The operators $D_q$ at a fixed
point $q$ form a finite-dimensional linear space, we denote it $\goth
D(q,M)$. The lemma~\mythelemma{7.4} provides the following result.
\mytheorem{7.3} Any differentiation $D$ of the algebra of extended tensor fields $\bold T(M)$ is represented as a field of differential operators
$D_q\in\goth D(q,M)$, one per each point $q\in N$. Conversely, each smooth 
field of tensorial first order differential operators is a differentiation
of the algebra $\bold T(M)$.
\endproclaim
The theorem~\mythetheorem{7.3} solves the problem of localization 
announced in the very beginning of this section. \pagebreak Relying 
on this theorem, one can formulate the following definition for a
differentiation of the algebra of extended tensor fields $\bold T(M)$.
\mydefinition{7.4} Let $N=T^{r_1\ldots\,r_Q}_{s_1\ldots\,s_Q}\!M$ 
be a composite tensor bundle over a smooth real manifold $M$. A
differentiation $D$ of the algebra $\bold T(M)$ is a smooth 
operator-valued function in $N$ such that it takes each point 
$q\in N$ to some differential operator $D_q\in\goth D(q,M)$.
\enddefinition
    Smoothness in both cases --- in theorem~\mythetheorem{7.3} and in
definition~\mythedefinition{7.4} means that the components $Z^i$, $Z^{h_1\ldots\,h_r}_{k_1\ldots\,k_s}[P]$, $\Gamma^k_i$ of the differential
operator $D_q$ are smooth functions of the variables \mythetag{3.3} in 
any local chart.\par
    Let's compare the definitions~\mythedefinition{7.4} and
\mythedefinition{4.1}. They are very similar. Therefore, by analogy to
the definition~\mythedefinition{1.3}, one can formulate the following
definition. 
\mydefinition{7.5}  A differentiation $D$ of the algebra $\bold T(M)$
is a geometric object in each local chart represented by its components
$Z^i$, $Z^{h_1\ldots\,h_r}_{k_1\ldots\,k_s}[P]$, $\Gamma^k_i$ and such 
that its components are smooth functions transformed according to the formulas \mythetag{7.14}, \mythetag{7.15}, \mythetag{7.16}, and 
\mythetag{7.17} under a change of local coordinates.
\enddefinition
     The definition~\mythedefinition{7.5} is a coordinate form 
of the definition~\mythedefinition{7.4}, and conversely, the definition~\mythedefinition{7.4} is a coordinate-free form of 
the definition~\mythedefinition{7.5}.
\head
8. Degenerate differentiations.
\endhead
\mydefinition{8.1} A differentiation $D$ of the algebra of extended 
tensor fields $\bold T(M)$ is called a {\it degenerate differentiation\/}
if its restriction $\delta\!:\,T^0_0(M)\to T^0_0(M)$ to the module $T^0_0(M)$ is identically zero.
\enddefinition
    For a degenerate differentiation the corresponding vector field 
\mythetag{7.2} is identically equal to zero. Then the equality
\mythetag{7.27} is obviously fulfilled, while the equality 
\mythetag{7.28} is reduced to the following one:
$$
\hskip -2em
\xi(\varphi\,\bold X)=\varphi\ \xi(\bold X).
\mytag{8.1}
$$
From the equality \mythetag{8.1} we conclude that the map $\xi\!:\,T^1_0(M)\to T^1_0(M)$ is an endomorphism of the 
module $T^1_0(M)$ over the ring $\goth F(N)$. 
\mydefinition{8.2} Let $A$ be a module over the ring of smooth 
real-valued functions $\goth F(M)$ in some manifold $M$. We say 
that the module $A$ admits a {\it localization\/} if it is 
isomorphic to a functional module so that each element $\bold a
\in A$ is represented as some function $\bold a(q)=\bold a_q$ in 
$M$ taking its values in some $\Bbb R$-linear spaces $A_q$ 
associated with each point $q$ of the manifold $M$.
\enddefinition
\mydefinition{8.3} Let $A$ be a module that admits a localization 
in the sense of the definition~\mythedefinition{8.2}. We say that
the localization of $A$ is a {\it complete localization} if the 
following two conditions are fulfilled:
\roster
\rosteritemwd=5pt
\item for any point $q\in M$ and for any vector $\bold v\in A_q$ 
      there exists an element $\bold a\in A$ such that $\bold a_q
      =\bold v$;
\item if $\bold a_q=0$ at some point $q\in M$, then there exist 
some finite set of elements $\bold E_0,\,\ldots,\,\bold E_n$ in
$A$ and some smooth functions $\alpha_0\,\,\dots,\,\alpha_n$ 
vanishing at the point $q$ such that $\bold a=\alpha_0\,\bold E_0
+\ldots+\alpha_n\,\bold E_n$.
\endroster      
\enddefinition
\mytheorem{8.1} Let $A$ and $B$ be two $\goth F(M)$-modules that
admit localizations. If the localization of $A$ is a complete
localization, then each homomorphism $f\!:\,A\to B$ is represented 
by a family of\/ $\Bbb R$-linear mappings 
$$
\hskip -2em
F_q\!:\,A_q\to B_q
\mytag{8.2}
$$
so that if\/ $\bold a\in A$ and $\bold b=f(\bold a)$, then $\bold b_q
=F_q(\bold a_q)$ for each point $q\in M$.
\endproclaim
\demo{Proof} First of all we should construct the mappings \mythetag{8.2}.
Let $q$ be some arbitrary point of the manifold $M$ and let $\bold v$ be
some arbitrary vector of the $\Bbb R$-linear space $A_q$. Then, according
to the item \therosteritem{1} of the definition~\mythedefinition{8.3}, we
have an element $\bold a\in A$ such that $\bold v=\bold a_q$. Applying the
homomorphism $f$ to $\bold a$ we get $\bold b=f(\bold a)\in B$. Let's 
define the mapping \mythetag{8.2} as follows:
$$
\hskip -2em
F_q(\bold v)=\bold b_q\text{, \ where \ }\bold b=f(\bold a)
\text{ \ for some \ }\bold a\in A\text{ \ such that \ }\bold v=\bold a_q.
\mytag{8.3}
$$
The choice of an element $\bold a\in A$ such  that $\bold v=\bold a_q$ is
not unique. Therefore, one should prove the consistence of the definition
\mythetag{8.3}. Suppose that $\bold a$ and $\tilde\bold a$ are two elements
of the module $A$ such that $\bold v=\bold a_q$ and $\bold v=\tilde\bold
a_q$. Then for the element $\bold c=\bold a-\tilde\bold a$ we get $\bold 
c_q=0$. Applying the item \therosteritem{2} of the 
definition~\mythedefinition{8.3} to $\bold c$, we get
$$
\hskip -2em
\bold c=\alpha_0\,\bold E_0+\ldots+\alpha_n\,\bold E_n,
\mytag{8.4}
$$
where $\alpha_0,\,\ldots,\,\alpha_n$ are smooth functions vanishing at the
point $q\in M$. Applying the homomorphism $f$ to \mythetag{8.4}, we derive
$$
\hskip -2em
\bold d=f(\bold a)-f(\tilde\bold a)=f(\bold c)=
\alpha_0\,f(\bold E_0)+\ldots+\alpha_n\,f(\bold E_n).
\mytag{8.5}
$$
Since $\alpha_0(q)=\ldots=\alpha_n(q)=0$, from \mythetag{8.5} we obtain
$\bold d_q=0$. This means that two elements $\bold b=f(\bold a)$ and 
$\tbb=f(\tilde\bold a)$ determine the same vector $\bold b_q=\tbb_q$
in \mythetag{8.3}. So, \mythetag{8.3} is a consistent definition of 
a mapping $F_q\!:\,A_q\to B_q$.\par
     The mapping $F_q\!:\,A_q\to B_q$ consistently defined by the 
equality \mythetag{8.3} is $\Bbb R$-linear. This fact is a trivial consequence of the equalities 
$$
\xalignat 2
&f(\bold a_1+\bold a_2)=f(\bold a_1)+f(\bold a_2)
&&\text{for any \ }\bold a_1,\bold a_2\in A,\\
&f(\lambda\,\bold a)=\lambda\,f(\bold a)
&&\text{for any \ }\bold a_2\in A\text{ \ and \ }
\lambda\in\goth F(M),
\endxalignat
$$
which mean that $f$ is a homomorphism of $\goth F(M)$-modules. And 
finally, from the equality \mythetag{8.3} it follows that for any
$\bold a\in A$ its image $\bold b=f(\bold a)$ is represented by a
function $\bold b$ whose values are obtained by applying the mappings
\mythetag{8.2} to the values of $\bold a$. The theorem is proved.
\qed\enddemo
\mytheorem{8.2} Let $\pi\!:\,VM\to M$ be a smooth $n$-dimensional vector
bundle over some base manifold $M$ and let $A$ be the set of all global
smooth sections\footnotemark\ of this bundle. Then $A$ admits complete
localization in the sense of the definition~\mythedefinition{8.3}.
\endproclaim
\footnotetext{ \ See the definition in \cite{31}.}
\adjustfootnotemark{-1}
\demo{Proof} The module structure of $A$ and its localization are obvious.
By definition, each section $\bold a$ of the bundle $\pi\!:\,VM\to M$ is a
function taking its values in fibers $V_q=\pi^{-1}(q)$ of the bundle $VM$.
We need to prove that this natural localization of $A$ is complete. Let's
begin with the item \therosteritem{1} in the
definition~\mythedefinition{8.3}. Let $q$ be some fixed point of the base
manifold $M$ and let $\bold v$ be some vector in the fiber $V_q$ over this
point. Each vector bundle is locally trivial. This means that there is some
open neighborhood $U$ of our point $q$ such that $\pi^{-1}(U)$ is isomorphic
to the trivial vector bundle $U\!\times\Bbb R^n$. This fact is expressed
by the following diagram:
$$
\hskip -2em
\CD
\pi^{-1}(U) @>\psi>> U\!\times\Bbb R^n\\
@V\pi VV @VVV\\
U @>\idop>> U.
\endCD
\mytag{8.6}
$$
The isomorphism $\psi$ in the diagram \mythetag{8.6} is linear within
each fiber of the bundle $VM$. Let's denote $\bold r_0=\psi(\bold v)\in
\Bbb R^n$ and let's choose the constant section $\bold r(q)\equiv\bold 
r_0$ of the trivial bundle $U\!\times\Bbb R^n$. Its preimage $\bold b=
\psi^{-1}(\bold r)$ is a smooth section of $VM$ over the open set $U$ 
and $\bold b_q=\bold v$. However, $\bold b$ is a local section. In order
to convert it to a global section let's choose some smooth function 
$\eta$ in $M$ such that $\eta(q)=1$ at the point $q$ and such that $\eta
\equiv 0$ outside the open set $U$. Then let's define 
$$
\hskip -2em
\bold a=\cases\eta\ \bold b &\text{within \ }U,\\
\ \ 0&\text{outside \ }U.\endcases
\mytag{8.7}
$$
It is easy to see that $\bold a$ in \mythetag{8.7} is a smooth global
section of the bundle $VM$. From $\eta(q)=1$ and froom $\bold b_q=\bold v$
we derive that $\bold a_q=\bold v$. Thus, we have proved that the module
$A$ fits the item \therosteritem{1} in the 
definition~\mythedefinition{8.3}.\par
     Now let's proceed with the item \therosteritem{2} in the 
definition~\mythedefinition{8.3}. Suppose that $\bold a$ is a smooth section
of the bundle $VM$ such that $\bold a_q=0$ at the point $q\in M$. Applying
the isomorphism $\psi$ taken from the diagram \mythetag{8.6} to the
restriction of $\bold a$ to the open set $U$, we get the following smooth
section of the trivial bundle $U\!\times\Bbb R^n$:
$$
\hskip -2em
\bold r=\psi(\bold a)=\Vmatrix \beta_1\\ \vdots\\ \beta_n\endVmatrix=
\sum^n_{i=1}\beta_i\,\bold e_i.
\mytag{8.8}
$$
Here $\bold e_1,\,\ldots,\,\bold e_n$ are constant unit vectors in 
$\Bbb R^n$:
$$
\hskip -2em
\bold e_1=\Vmatrix 1\\ 0\\ \vdots\\ 0\\ 0\endVmatrix,\quad
\bold e_2=\Vmatrix 0\\ 1\\ \vdots\\ 0\\ 0\endVmatrix,\quad
.\ .\ .\ ,\ \quad
\bold e_{n-1}=\Vmatrix 0\\ 0\\ \vdots\\ 1\\ 0\endVmatrix,\quad
\bold e_n=\Vmatrix 0\\ 0\\ \vdots\\ 0\\ 1\endVmatrix.
\mytag{8.9}
$$
From $\bold a_q=0$ and from \mythetag{8.8} we derive that the smooth
functions $\beta_1,\,\ldots,\,\beta_n$ should vanish at the point $q$,
i\.\,e\. $\beta_i(q)=0$. In \mythetag{8.8} the vectors \mythetag{8.9}
represent some constant smooth sections of the trivial bundle $U\!\times
\Bbb R^n$. Let's denote $\hat\bold E_i=\psi^{-1}(\bold e_i)$. Then $\hat
\bold E_1,\,\ldots,\,\hat\bold E_n$ are smooth local sections of the 
bundle $VM$ over the open set $U$. From \mythetag{8.8} we derive the
following expansion:
$$
\hskip -2em
\bold a=\sum^n_{i=1}\beta_i\,\hat\bold E_i.
\mytag{8.10}
$$
The expansion \mythetag{8.10} is a local expansion, it is not the
expansion of $\bold a $, but of its restriction to the open set $U$.
In order to make it global let's define
$$
\xalignat 2
&\hskip -2em
\alpha_i=\cases\eta\ \beta_i &\text{within \ }U,\\
\ \ 0&\text{outside \ }U,\endcases
&&\bold E_i=\cases\eta\ \hat\bold E_i &\text{within \ }U,\\
\ \ 0&\text{outside \ }U.\endcases\quad
\mytag{8.11}
\endxalignat
$$
Here $\eta$ again is a smooth function in $M$ such that $\eta(q)=1$
and such that it is identically zero outside $U$. Multiplying both 
sides of the expansion \mythetag{8.10} by $\eta^2$ and taking into 
account \mythetag{8.11}, we derive
$$
\hskip -2em
\bold a\ \eta^2=\sum^n_{i=1}\alpha_i\,\bold E_i.
\mytag{8.12}
$$
Then, for the original section $\bold a$ due to 
\mythetag{8.12} we get
$$
\bold a=(1-\eta^2)\,\bold a+\bold a\ \eta^2=
(1-\eta^2)\,\bold a+\sum^n_{i=1}\alpha_i\,\bold E_i.
$$
Let's denote $\alpha_0=1-\eta$ and $\bold E_0=(1+\eta)\,\bold a$. Then
the above local expansion \mythetag{8.10} is transformed to the following
global one:
$$
\hskip -2em
\bold a=\sum^n_{i=0}\alpha_i\,\bold E_i.
\mytag{8.13}
$$
Since $\beta_i(q)=0$ for $i=1,\,\ldots,\,n$ and since $\eta(q)=1$, we find
that all coefficients $\alpha_i$ in \mythetag{8.13} do vanish at the point
$q$. Comparing \mythetag{8.13} with the expansion in the item
\therosteritem{2} of the definition~\mythedefinition{8.3}, we finally
conclude that the module $A$ admits a complete localization. The
theorem~\mythetheorem{8.2} is proved.
\qed\enddemo
    Now let's return to the definition~\mythedefinition{8.1} and consider
the map $\xi\!:\,T^1_0(M)\to T^1_0(M)$ associated with some differentiation
$D$. The extended vector fields are naturally interpreted as sections of
a vector bundle. Indeed, any tensor bundle is a vector bundle in the sense
that any tensor space $T^\alpha_\beta(p,M)$ is a linear vector space over the real numbers $\Bbb R$. Therefore, one can apply the 
theorem~\mythetheorem{8.2} to the module $T^1_0(M)$. As it was mentioned above, the mapping $\xi\!:\,T^1_0(M)\to T^1_0(M)$ is an endomorphism.
Applying the theorem~\mythetheorem{8.1} to it, we find that $\xi$ is given
by some extended tensor field $\bold S$ of the type $(1,1)$ acting as a
linear operator at each point $q\in N$:	
$$
\hskip -2em
\xi(\bold X)=C(\bold S\otimes\bold X)\text{ \ for any \ }
\bold X\in T^1_0(M).
\mytag{8.14}
$$
\mytheorem{8.3} Defining a degenerate differentiation $D$ of the 
algebra $\bold T(M)$ is equivalent to defining some extended tensor 
field\/ $\bold S$ of the type $(1,1)$.
\endproclaim
     Apart from \mythetag{8.14}, the theorem~\mythetheorem{8.3} 
can be understood in a coordinate form. Indeed, for a degenerate
differentiation $D$ from \mythetag{7.2} and from the 
definition~\mythedefinition{8.1} we derive $Z^i=0$ and
$Z^{i_1\ldots\,i_r}_{j_1\ldots\,j_s}[P]=0$. Since $Z^i=0$, the 
transformation formulas \mythetag{7.16} and \mythetag{7.17} now
are written as follows:
$$
\pagebreak
\xalignat 2
&\hskip -2em
\Gamma^k_i=\sum^n_{b=1}\sum^n_{a=1}S^k_a\,T^b_i\ \tilde\Gamma^a_b,
&&\tilde\Gamma^k_i=\sum^n_{b=1}\sum^n_{a=1}T^k_a\,S^b_i\ \Gamma^a_b.
\qquad
\mytag{8.15}
\endxalignat
$$
The formulas \mythetag{8.15} mean that $\Gamma^k_i$ and $\tilde
\Gamma^a_b$ are the components of some extended tensor field of
the type $(1,1)$, see \mythetag{1.12} and \mythetag{1.13} for
comparison. Applying \mythetag{7.26} to \mythetag{8.15}, we find
that $\Gamma^k_i$ are the components of the tensor field $\bold S$ 
in some local chart. This is a coordinate proof for the 
theorem~\mythetheorem{8.3}.\par
\head
9. Covariant differentiations.
\endhead
     The set of differentiations of the extended algebra $\bold T(M)$
possesses the structure of a module over the ring of smooth functions
$\goth F(N)$. The set of extended vector fields $T^1_0(M)$ is also a
module over the same ring $\goth F(N)$. Therefore, the following 
definition is consistent.
\mydefinition{9.1} Say that in a manifold $M$ a covariant differentiation
of the algebra of extended tensor fields $\bold T(M)$ is given if some
homomorphism of $\goth F(N)$-modules \linebreak $\nabla\!:\,T^1_0(M)\to\goth D(M)$ is given. The image of a vector field
$\bold Y$ under such homomorphism is denoted by $\nabla_{\bold Y}$. The
differentiation $D=\nabla_{\bold Y}\in\goth D(M)$ is called the {\it
covariant differentiation along the vector field\/} $\bold Y$.
\enddefinition
     Let's remember that the module $\goth D(M)$ admits a localization.
Indeed, according to the theorem~\mythetheorem{7.3} each differentiation
$D$ is a field of differential operators. As for the differential operators
themselves, they form finite-dimensional $\Bbb R$-linear spaces $\goth 
D(q,M)$, one per each point $q\in N$. The coefficients $Z^i$, $Z^{i_1\ldots
\,i_r}_{j_1\ldots\,j_s}[P]$, and $\Gamma^k_i$ from \mythetag{7.26} are coordinates within $\goth D(q,M)$. Therefore, we have
$$
\hskip -2em
\dim(\goth D(q,M))=n+n^2+\sum^Q_{P=1}n^{r_P+s_P}=n^2+\dim N.
\mytag{9.1}
$$
Under a change of a local chart the coefficients $Z^i$, $Z^{i_1\ldots
\,i_r}_{j_1\ldots\,j_s}[P]$, $\Gamma^k_i$ are transformed according to
the formulas \mythetag{7.14}, \mythetag{7.15}, \mythetag{7.16}, and 
\mythetag{7.17}. These formulas are linear with respect to $Z^i$, $Z^{i_1
\ldots\,i_r}_{j_1\ldots\,j_s}[P]$, $\Gamma^k_i$ and with respect to
transformed coordinates $\tilde Z^i$, $\tilde Z^{i_1\ldots\,i_r}_{j_1
\ldots\,j_s}[P]$, $\tilde \Gamma^k_i$. Therefore, the spaces $\goth D(q,M)$
are glued into a vector bundle of the dimension \mythetag{9.1} for which
$N$ is a base manifold. This fact means that one can apply the theorem~\mythetheorem{8.2} to the module $\goth D(M)$.\par
    Let $\nabla$ be some covariant differentiation of the algebra of
extended tensor fields. Then, applying the theorem~\mythetheorem{8.1} to
the homomorphism $\nabla\!:\,T^1_0(M)\to\goth D(M)$, we find that this
homomorphism is composed by $R$-linear maps $T_{\pi(q)}(M)\to\goth D(q,M)$
specific to each point $q\in N$. This fact is expressed by the following
formula:
$$
\hskip -2em
\nabla_{\bold Y}\bold X=C(\bold Y\otimes\nabla\bold X)
=\sum^n_{j=1}
\dsize\msum{n}\Sb i_1,\,\ldots,\,i_\alpha\\j_1,\,\ldots,\,j_\beta\endSb
Y^j\ \nabla_{\!\!j}X^{i_1\ldots\,i_\alpha}_{j_1\ldots\,j_\beta}\ 
\bold E^{j_1\ldots\,j_\beta}_{i_1\ldots\,i_\alpha}.
\mytag{9.2}
$$
Looking at \mythetag{9.2}, we see that each covariant differentiation
$\nabla$ can be treated as an operator producing the extended tensor 
field $\nabla\bold X$ of the type $(\alpha,\beta+1)$ from any given
extended tensor field $\bold X$ of the type $(\alpha,\beta)$. This
operator increases by one the number of covariant indices of a tensor
field $\bold X$. It is called the operator of {\it covariant 
differential\/} associated with the covariant differentiation 
$\nabla$.\par
     Let's consider the linear map $T_{\pi(q)}(M)\to\goth D(q,M)$ 
produced by some covariant differentiation $\nabla$ at some particular
point $q\in N$. In a local chart this map is given by some linear functions
expressing the components of the differential operator $D_q$, where
$D=\nabla_{\bold Y}$, through the components of the vector $\bold Y_{\!q}$:
$$
\align
&\hskip -2em
Z^{i_1\ldots\,i_r}_{j_1\ldots\,j_s}[P]=\sum^n_{j=1}
Z^{\ i_1\ldots\,i_r}_{j\,j_1\ldots\,j_s}[P]\ Y^j,
\mytag{9.3}\\
&\hskip -2em
Z^i=\sum^n_{j=1}Z^i_j\ Y^j,
\qquad
\Gamma^k_i=\sum^n_{j=1}\Gamma^k_{j\,i}\ Y^j.
\mytag{9.4}
\endalign
$$
Substituting \mythetag{9.3} and \mythetag{9.4} into the formula \mythetag{7.26}, then taking into account \mythetag{9.2} 
and \mythetag{7.25}, we derive the following formula:
$$
\hskip 1pt
\gathered
\nabla_{\!\!j}X^{i_1\ldots\,i_\alpha}_{j_1\ldots\,j_\beta}
=\sum^n_{i=1}Z^i_j\ \frac{\partial X^{i_1\ldots\,i_\alpha}_{j_1
\ldots\,j_\beta}}{\partial x^i}+\sum^Q_{P=1}\dsize\msum{n}\Sb h_1,\,\ldots,\,h_r\\k_1,\,\ldots,\,k_s\endSb Z^{\,h_1\ldots\,h_r}_{j\,k_1\ldots\,k_s}[P]
\ \frac{\partial X^{i_1\ldots\,i_\alpha}_{j_1\ldots\,j_\beta}}
{\partial T^{h_1\ldots\,h_r}_{k_1\ldots\,k_s}[P]
\vphantom{\vrule height 10pt depth 0pt}}\ +\\
+\sum^\alpha_{m=1}\sum^n_{v_m=1}\Gamma^{i_m}_{j\,v_m}\ X^{\,i_1\ldots\,v_m
\ldots\,i_\alpha}_{j_1\ldots\,\ldots\,\ldots\,j_\beta}-\sum^\beta_{m=1}
\sum^n_{w_m=1}\Gamma^{w_m}_{j\,j_m}\ X^{\,i_1\ldots\,\ldots\,\ldots\,
i_\alpha}_{j_1\ldots\,w_m\ldots\,j_\beta},
\endgathered
\mytag{9.5}
$$
From \mythetag{7.14} and \mythetag{7.15} we derive the following transformation formulas for the quantities $Z^i_j$, $Z^{\ i_1\ldots
\,i_r}_{j\,j_1\ldots\,j_s}[P]$ in \mythetag{9.3} and \mythetag{9.4}:
$$
\align
&\cases
Z^i_j=\dsize\sum^n_{h=1}\sum^n_{k=1}S^i_h\,T^k_j\ \tilde Z^h_k,\\
\vspace{2ex}
\gathered
Z^{\,i_1\ldots\,i_r}_{j\,j_1\ldots\,j_s}[P]=\sum^n_{k=1}
\dsize\msum{n}\Sb h_1,\,\ldots,\,h_r\\k_1,\,\ldots,\,k_s\endSb 
S^{i_1}_{h_1}\ldots\,S^{i_r}_{h_r}
\ T^{k_1}_{j_1}\ldots\,T^{k_s}_{j_s}\,T^k_j\ 
\tilde Z^{\,h_1\ldots\,h_r}_{k\,k_1\ldots\,k_s}[P]\ -\\
-\sum^r_{m=1}\sum^n_{i=1}\sum^n_{h=1}\sum^n_{k=1}\sum^n_{v_m=1}
\theta^{i_m}_{i\,v_m}\,T^{\,i_1\ldots\,v_m\ldots\,i_r}_{j_1
\ldots\,\ldots\,\ldots\,j_s}[P]\ S^i_h\,T^k_j\ \tilde Z^h_k\ +\\
+\sum^s_{m=1}\sum^n_{i=1}\sum^n_{h=1}\sum^n_{k=1}\sum^n_{w_m=1}
\theta^{w_m}_{i\,j_m}\,T^{\,i_1\ldots\,\ldots\,\ldots\,
i_r}_{j_1\ldots\,w_m\ldots\,j_s}[P]\ S^i_h\,T^k_j\ \tilde Z^h_k,
\endgathered
\endcases\ \qquad
\mytag{9.6}\\
\vspace{2ex}
&\cases
\tilde Z^i_j=\dsize\sum^n_{h=1}\sum^n_{k=1}T^i_h\,S^k_j\ Z^h_k,\\
\vspace{2ex}
\gathered
\tilde Z^{\,i_1\ldots\,i_r}_{j\,j_1\ldots\,j_s}[P]=\sum^n_{k=1}
\dsize\msum{n}\Sb h_1,\,\ldots,\,h_r\\k_1,\,\ldots,\,k_s\endSb 
T^{i_1}_{h_1}\ldots\,T^{i_r}_{h_r}
\ S^{k_1}_{j_1}\ldots\,S^{k_s}_{j_s}\,S^k_j\ 
Z^{\,h_1\ldots\,h_r}_{k\,k_1\ldots\,k_s}[P]\ -\\
-\sum^r_{m=1}\sum^n_{i=1}\sum^n_{h=1}\sum^n_{k=1}\sum^n_{v_m=1}
\tilde\theta^{i_m}_{i\,v_m}\,T^{\,i_1\ldots\,v_m\ldots\,i_r}_{j_1
\ldots\,\ldots\,\ldots\,j_s}[P]\ T^i_h\,S^k_j\ Z^h_k\ +\\
+\sum^s_{m=1}\sum^n_{i=1}\sum^n_{h=1}\sum^n_{k=1}\sum^n_{w_m=1}
\tilde\theta^{w_m}_{i\,j_m}\,T^{\,i_1\ldots\,\ldots\,\ldots\,
i_r}_{j_1\ldots\,w_m\ldots\,j_s}[P]\ T^i_h\,S^k_j\ Z^h_k.
\endgathered
\endcases\ \qquad
\mytag{9.7}
\endalign
$$
Similarly, from \mythetag{7.16} and \mythetag{7.17} we derive the
transformation formulas for the quantities $\Gamma^k_{j\,i}$ in
the formula \mythetag{9.4}:
$$
\align
&\hskip -2em
\Gamma^k_{j\,i}=\dsize\sum^n_{b=1}\sum^n_{a=1}\sum^n_{c=1}
S^k_a\,T^b_i\,T^c_j\ \tilde\Gamma^a_{c\,b}+
\sum^n_{a=1}Z^a_j\,\theta^k_{ai},
\mytag{9.8}\\
\vspace{2ex}
&\hskip -2em
\tilde\Gamma^k_{j\,i}=\dsize\sum^n_{b=1}\sum^n_{a=1}\sum^n_{c=1}
T^k_a\,S^b_i\,S^c_j\ \Gamma^a_{c\,b}+
\sum^n_{a=1}\tilde Z^a_j\,\tilde\theta^k_{ai}.
\mytag{9.9}
\endalign
$$
The formula \mythetag{9.5} yields the explicit expression for an
arbitrary covariant derivative in general case. However, below
we consider some specializations of this formula which appear to
be more valuable than the formula \mythetag{9.5} itself.
\head
10. Degenerate covariant differentiations.
\endhead
\mydefinition{10.1} A covariant differentiation $\nabla$ is said to 
be {\it degenerate\/} if \ $\nabla_{\bold Y}\psi=0$ for any extended
scalar field $\psi$ and for any extended vector field $\bold Y$.
\enddefinition
     This definition is concordant with the
definition~\mythedefinition{8.1}. If\/ $\nabla$ is a degenerate
covariant differentiation, then $D=\nabla_{\bold Y}$ is a degenerate
differentiation for any extended vector field $\bold Y$. According to
the theorem~\mythetheorem{8.3} and formula \mythetag{8.14}, $D$ is
associated with some extended tensor field $\bold S$ of the type
$(1,1)$.  In the present case this field should depend of $\bold Y$,
so the homomorphism $\nabla\!:\,T^1_0(M)\to\goth D(M)$ in the 
definition~\mythedefinition{9.1} reduces to the homomorphism
$$
\hskip -2em
T^1_0(M)\to T^1_1(M).
\mytag{10.1}
$$
Applying the theorem~\mythetheorem{8.1} to the homomorphism 
\mythetag{10.1} we derive the following theorem for degenerate
covariant differentiations.
\mytheorem{10.1} Defining a degenerate covariant differentiation
$\nabla$ of the algebra $\bold T(M)$ is equivalent to defining some
extended tensor field\/ $\bold S$ of the type $(1,2)$.
\endproclaim
     Like the theorem~\mythetheorem{8.3}, this theorem can be 
understood in a coordinate form. Indeed, if $\nabla$ is degenerate, 
the vector field \mythetag{7.2} associated with the differentiation
$D=\nabla_{\bold Y}$ should be identically zero for any extended
vector field $\bold Y$. This means that the coefficients 
$Z^{\ i_1\ldots\,i_r}_{j\,j_1\ldots\,j_s}[P]$ and $Z^i_j$ in
\mythetag{9.3} and \mythetag{9.4} are equal to zero. Then the
transformation formulas \mythetag{9.8} and \mythetag{9.9} are
written as follows:
$$
\align
&\hskip -2em
\Gamma^k_{j\,i}=\dsize\sum^n_{b=1}\sum^n_{a=1}\sum^n_{c=1}
S^k_a\,T^b_i\,T^c_j\ \tilde\Gamma^a_{c\,b},
\mytag{10.2}\\
\vspace{2ex}
&\hskip -2em
\tilde\Gamma^k_{j\,i}=\dsize\sum^n_{b=1}\sum^n_{a=1}\sum^n_{c=1}
T^k_a\,S^b_i\,S^c_j\ \Gamma^a_{c\,b}.
\mytag{10.3}
\endalign
$$
Comparing \mythetag{10.2} and \mythetag{10.3} with \mythetag{1.12} and
\mythetag{1.13}, we see that $\Gamma^k_{j\,i}$ behave like the components
of a tensor of the type $(1,2)$. They are that very quantities that define
the extended tensor field $\bold S$ in a local chart.
\head
11. Horizontal and vertical covariant differentiations.
\endhead
     Suppose again that some composite tensor bundle $N=T^{r_1\ldots
\,r_Q}_{s_1\ldots\,s_Q}\!M$ over a base manifold $M$ is fixed. Let
$\nabla$ be a covariant differentiation of the algebra of extended
tensor fields $\bold T(M)$. Then $D=\nabla_{\bold Y}$ is a differentiation
of $\bold T(M)$, its restriction to the set of scalar fields is given by
some vector field $\bold Z=\bold Z(\bold Y)$ in $N$ (see formula \mythetag{7.2} above). In other words, we have a homomorphism
$$
\hskip -2em
T^1_0(M)\to T^1_0(N)
\mytag{11.1}
$$
that maps an extended vector field $\bold Y$ of $M$ to some regular
vector field of $N$. Applying the localization theorem~\mythetheorem{8.1}
to the homomorphism \mythetag{11.1}, we come to the following definition
and the theorem after it.
\mydefinition{11.1} Suppose that for each point $q$ of the composite 
tensor bundle $N$ over the base $M$ some $\Bbb R$-linear map of the 
vector spaces 
$$
\hskip -2em
f\!:\ T_{\pi(q)}(M)\to T_q(N)
\mytag{11.2}
$$ is given. Then we say that a {\it lift\/} of vectors from $M$ to the bundle $N$ is defined.
\enddefinition
\mytheorem{11.1} Any homomorphism of $\goth F(N)$-modules \mythetag{11.1}
is uniquely associated with some smooth lift of vectors from $M$ to $N$.
It is represented by this lift as a collection of $\Bbb R$-linear maps
$T_{\pi(q)}(M)\to T_q(N)$ specific to each point $q\in N$.
\endproclaim
    Now let's consider the canonical projection $\pi\!:\,N\to M$. The
differential of this map acts in the direction opposite to the lift of vectors \mythetag{11.2}  introduced in the
definition~\mythedefinition{11.1}. Indeed, we have $\pi_*\!:\,T_q(N)\to
T_{\pi(q)}(M)$ at each point $q\in N$. Therefore, the composition
$f\compos\pi_*$ acts from $T_{\pi(q)}(M)$ to $T_{\pi(q)}(M)$. This
composite map determines an extended operator field (a tensor field 
of the type $(1,1)$).
\mydefinition{11.2} A lift of vectors $f$ from $M$ to $N$ is called
{\it vertical\/} if $\pi_*\compos f=0.$
\enddefinition
\mydefinition{11.3} A lift of vectors $f$ from $M$ to $N$ is called
{\it horizontal\/} if $\pi_*\compos f=\idop$, i\.\,e\. if the composition
$\pi_*\compos f$ coincides with the field of identical operators.
\enddefinition
     Like any other bundle, the composite tensor bundle $N$ naturally
subdivides into fibers over the points of the base manifold $M$.
The set of vectors tangent to the fiber at a point $q$ is a linear subspace
within the tangent space $T_q(N)$. This subspace coincides with the kernel
of the mapping $\pi_*$. We denote this subspace 
$$
\hskip -2em
V_q(N)=\Ker\pi_*\subset T_q(N)
\mytag{11.3}
$$ 
and call it the {\it vertical subspace}. Any vertical lift of vectors determines a set linear maps from $T_{\pi(q)}$ to the vertical subspace 
\mythetag{11.3} for each point $q\in N$.
\mylemma{11.1} The difference of two horizontal lifts is a vertical lift 
of vectors from the base manifold $M$ to the bundle $N$.
\endproclaim
Indeed, if one takes two horizontal lifts of vectors $f_1$ and $f_2$, 
then $\pi_*\compos(f_1-f_2)=\ =\pi_*\compos f_1-\pi_*\compos f_2=\idop
-\idop=0$. This means that the difference $f_1-f_2$ is a vertical lift
according to the definition~\mythedefinition{11.2}.\par
    Each covariant differentiation $\nabla$ is associated with some
lift of vectors (see the definition~\mythedefinition{11.1}, and the theorem~\mythetheorem{11.1} above).
\mydefinition{11.4} A covariant differentiation $\nabla$ of
the algebra of extended tensor fields $\bold T(M)$ is called a {\it horizontal differentiation\/} (or a {\it vertical differentiation\/}) 
if the corresponding lift of vectors is {\it horizontal\/} (or {\it
vertical}).
\enddefinition
\mylemma{11.2} The difference of two horizontal covariant differentiations
is a vertical covariant differentiation of the algebra of extended tensor
fields $\bold T(M)$.
\endproclaim
    The lemma~\mythelemma{11.2} is an immediate consequence of the
lemma~\mythelemma{11.1}. Unlike \mycite{4}, here we shall not pay much
attention to vertical covariant differentiations. In the present more
general theory they are replaced by a more general construct.
\head
12. Native extended tensor fields\\
and vertical multivariate differentiations.
\endhead
     Let $N=T^{r_1\ldots\,r_Q}_{s_1\ldots\,s_Q}\!M$ be a composite tensor
bundle over a base manifold $M$. Then each its point $q$ is represented 
by a list $q=(p,\,\bold T[1],\,\ldots,\,\bold T[Q])$, where $p\in M$ and
$\bold T[1],\,\ldots,\,\bold T[Q]$ are some tensors at the point $p$ (see
formula \mythetag{3.2} above). Let's consider the map that takes $q$ to
the $P$-th tensor $\bold T[P]$ in the list. According to the 
definition~\mythedefinition{4.1}, this map is an extended tensor field of
the type $(r_P,s_P)$. It is canonically associated with the bundle $N$.
Therefore, it is called a {\it native extended tensor field}. Totally, we
have $Q$ native extended tensor fields associated with the composite
tensor bundle $N=T^{r_1\ldots\,r_Q}_{s_1\ldots\,s_Q}\!M$, we denote 
them $\bold T[1],\,\ldots,\,\bold T[Q]$.\par
\mydefinition{12.1} A {\it multivariate differentiation\/} of the type
$(\beta,\alpha)$ in the algebra of extended tensor fields $\bold T(M)$ 
is a homomorphism of $\goth F(N)$-modules
$$
\hskip -2em
\nabla\!:\,T^\alpha_\beta(M)\to\goth D(M).
\mytag{12.1}
$$
If\/ $\bold Y$ is an extended tensor field of the type $(\alpha,\beta)$,
then we can apply the homomorphism \mythetag{12.1} to it. As a result 
we get the differentiation $D=\nabla_{\bold Y}$ of the algebra
$\bold T(M)$. It is called the {\it multivariate differentiation along
the tensor field\/} $\bold Y$.
\enddefinition
     Note that the type of a multivariate differentiation $(\beta,\alpha)$
in the above definition~\mythedefinition{12.1} is dual to the type of the
module $T^\alpha_\beta(M)$ in the formula \mythetag{12.1}. If $\alpha=1$ 
and $\beta=0$, the definition~\mythedefinition{12.1} reduces to the 
definition~\mythedefinition{9.1}. This means that a covariant
differentiation is a special multivariate differentiation whose type
is $(0,1)$. Similarly, a multivariate differentiation of the type $(1,0)$
is called a {\it contravariant differentiation}. Covariant differentiations
of extended tensor fields appear to be a useful tool in describing Newtonian
dynamical systems in Riemannian manifolds (see \myciterange{4}{4--23}). 
The same is true for contravariant differentiations in the case of 
Hamiltonian dynamical systems (see \myciterange{24}{24--30}). As for 
general multivariate differentiations introduced in the above
definition~\mythedefinition{11.1}, I think they will find their proper
place in theories of continuous media (see \mycite{34} and 
\myciterange{36}{36--39}) and in field theories.\par
     {\bf A remark}. Let's consider the special case, where the tensor field
$\bold Y$ of the type $(\alpha,\beta)$ is constructed as a tensor product:
$$
\hskip -2em
\bold Y=\bold Y[1]\otimes\ldots\otimes \bold Y[\alpha]
\otimes\bold H[1]\otimes\ldots\otimes\bold H[\beta].
\mytag{12.2}
$$
Here $\bold Y[1],\,\ldots,\,\bold Y[\alpha]$ are some vector fields and
$\bold H[1],\,\ldots,\,\bold H[\beta]$ are some covector fields.
Substituting \mythetag{12.2} into $\nabla_{\bold Y}$, we find that 
$\nabla_{\bold Y}$ is a differentiation depending on $\alpha$ vectorial
variables and $\beta$ covectorial variables. Keeping in mind this special
case, we used the term {\tencyr\char '074}multivariate 
differentiation{\tencyr\char '076} for $\nabla_{\bold Y}$ in the
definition~\mythedefinition{12.1}.\par
    Let $\nabla$ be some multivariate differentiation of the algebra 
of extended tensor fields $\bold T(M)$. Then, applying the localization
theorem~\mythetheorem{8.1} to the homomorphism $\nabla\!:\,T^\alpha_\beta(M)
\to\goth D(M)$, we find that this homomorphism is composed by 
$\Bbb R$-linear maps $T^\alpha_\beta(\pi(q),M)\to\goth D(q,M)$ specific to
each point $q\in N$. This fact is expressed by the following formula similar
to the formula \mythetag{9.2} above:
$$
\hskip -2em
\gathered
\bold Y\mapsto\nabla_{\bold Y}\bold X=
\dsize\msum{n}\Sb i_1,\,\ldots,\,i_\alpha\\j_1,\,\ldots,\,j_\beta\endSb
\dsize\msum{n}\Sb h_1,\,\ldots,\,h_r\\k_1,\,\ldots,\,k_s\endSb
Y^{h_1\ldots\,h_r}_{k_1\ldots\,k_s}
\ \nabla^{k_1\ldots\,k_s}_{\!h_1\ldots\,h_r}\!X^{i_1\ldots\,i_\alpha}_{j_1
\ldots\,j_\beta}\ \times\\
\times\ \bold E^{j_1\ldots\,j_\beta}_{i_1\ldots\,i_\alpha}=
C(\bold Y\otimes\nabla\bold X).
\endgathered
\mytag{12.3}
$$
Looking at \mythetag{12.3}, we see that each multivariate differentiation
$\nabla$ of the type $(s,r)$ can be treated as an operator producing the
extended tensor field $\nabla\bold X$ of the type $(\alpha+s,\beta+r)$ 
from any given extended tensor field $\bold X$ of the type $(\alpha,
\beta)$. This operator is called the operator of {\it multivariate
differential\/} of the type $(s,r)$.\par
     Let $P$ be an integer number such that $1\leqslant P\leqslant Q$ 
and let $\bold Y$ be an extended tensor field of the type $(r_P,s_P)$.
Remember that each point $q$ of the composite tensor bundle $N=T^{r_1
\ldots\,r_Q}_{s_1\ldots\,s_Q}\!M$ is a list of the form \mythetag{3.2}:
$$
\hskip -2em
q=(p,\,\bold T[1],\,\ldots,\,\bold T[Q]).
\mytag{12.4}
$$
Note that the $P$-th tensor $\bold T[P]$ in the list \mythetag{12.4} 
has the same type as the tensor $\bold Y=\bold Y_{\!q}$ (the value of 
the extended tensor field $\bold Y$ at the point $q$). They both belong 
to the same tensor space $T^{r_P}_{s_P}(p,M)$, therefore we can add them.
Then
$$
\hskip -2em
q(t)=(p,\,\bold T[1],\,\ldots,\,\bold T[P]+t\,\bold Y_{\!q},
\,\ldots,\,\bold T[Q])
\mytag{12.5}
$$
is a one-parametric set of points in $N$, the scalar variable $t$ being 
its parameter. In other words, in \mythetag{12.5} we have a line (a 
straight line) passing through the initial point $q\in N$ and lying completely within the fiber over the point $p=\pi(q)\in M$. Suppose that
$\bold X$ is some extended tensor field of the type $(\alpha,\beta)$. Denote by $\bold X(t)$ the values of this tensor field at the points of the above parametric line \mythetag{12.5}:
$$
\hskip -2em
\bold X(t)=\bold X_{q(t)}.
\mytag{12.6}
$$
Since $\pi(q(t))=p=\const$ for any $t$, the values of the tensor-valued
function \mythetag{12.6} all belong to the same tensor space
$T^\alpha_\beta(p,M)$. Hence, we can add and subtract them, and, since $\bold X$ is a smooth field, we can take the following limit of the ratio:
$$
\hskip -2em
\dot X(t)=\lim_{\tau\to\,0}\frac{\bold X(t+\tau)-\bold X(t)}{\tau}.
\mytag{12.7}
$$
Let's denote by $\bold Z_q$ the value of the derivative \mythetag{12.7} 
for $t=0$:
$$
\hskip -2em
\bold Z_q=\dot X(0)=\frac{d\bold X_{q(t)}}{dt}\,\hbox{\vrule
height14pt depth8pt width 0.5pt}_{\ t=0}.\hskip -2em
\mytag{12.8}
$$
It is easy to understand that, \pagebreak when $q$ is fixed, $\bold Z_q$
is a tensor from the tensor space $T^\alpha_\beta(p,M)$ at the point $p=\pi(q)$. By varying $q\in N$, we find that the tensors $\bold Z_q$
constitute a smooth extended tensor field $\bold Z$. So, we have constructed
a map
$$
\hskip -2em
D\!:\,\bold T(M)\to\bold T(M).
\mytag{12.9}
$$
It is easy to check up that the map \mythetag{12.9} defined by means 
of the formulas \mythetag{12.5}, \mythetag{12.6}, \mythetag{12.7}, and \mythetag{12.8} is a differentiation of the algebra of extended tensor
fields $\bold T(M)$, i\.\,e\. $D\in\goth D(M)$ (see definition~\mythedefinition{6.2} above). Moreover, due to the formula
\mythetag{12.5} this differentiation $D$ depends on the extended tensor
field $\bold Y$. This dependence $D=D(\bold Y)$ is a homomorphism $T^{r_P}_{s_P}(M)\to\goth D(M)$ fitting the
definition~\mythedefinition{12.1}. The easiest way to prove this fact is
to write the equality \mythetag{12.8} in a local chart, i\.\,e\. in some
local coordinates \mythetag{3.3}:
$$
\hskip -2em
\bold Z=
\dsize\msum{n}\Sb i_1,\,\ldots,\,i_\alpha\\j_1,\,\ldots,\,j_\beta\endSb
\dsize\msum{n}\Sb h_1,\,\ldots,\,h_r\\k_1,\,\ldots,\,k_s\endSb
Y^{h_1\ldots\,h_r}_{k_1\ldots\,k_s}\
\frac{\partial X^{i_1\ldots\,i_\alpha}_{j_1\ldots\,j_\beta}}
{\partial T^{h_1\ldots\,h_r}_{k_1\ldots\,k_s}[P]
\vphantom{\vrule height 10pt depth 0pt}}
\ \bold E^{j_1\ldots\,j_\beta}_{i_1\ldots\,i_\alpha}.
\mytag{12.10}
$$
Here $r=r_P$ and $s=s_P$. Now, comparing \mythetag{12.10} with the formula
\mythetag{12.3}, we can write $\bold Z=\vnabla_{\!\bold Y}[P]\bold X$, where
$\vnabla$ is a special sign, the {\tencyr\char '074}double bar 
nabla{\tencyr\char '076}, that we shall use for denoting the multivariate
differentiations defined through the formulas \mythetag{12.5},
\mythetag{12.6}, \mythetag{12.7}, and \mythetag{12.8}. In a local chart
$\vnabla[P]$ is represented by the formula
$$
\hskip -2em
\vnabla^{k_1\ldots\,k_s}_{h_1\ldots\,h_r}[P]=\frac{\partial }
{\partial T^{h_1\ldots\,h_r}_{k_1\ldots\,k_s}[P]
\vphantom{\vrule height 10pt depth 0pt}},
\mytag{12.11}
$$
where $r=r_P$ and $s=s_P$. This formula \mythetag{12.11} is a short
version of the formula \mythetag{12.10}. Following the tradition, we
shall use the term {\it multivariate derivative} for the differential
operator representing the differentiation $\vnabla[P]$ in local
coordinates.
\mydefinition{12.2} The multivariate differentiation $\vnabla[P]$ defined
through the formulas \mythetag{12.5}, \mythetag{12.6}, \mythetag{12.7},
\mythetag{12.8} and represented by the formula \mythetag{12.11} in local coordinates is called the {\it $P$-th canonical\footnotemark\ vertical multivariate differentiation} associated with the composite tensor bundle
$N=T^{r_1\ldots\,r_Q}_{s_1\ldots\,s_Q}\!M$.
\enddefinition
\footnotetext{ \ Note that $\vnabla[P]$ is canonically associated with 
the bundle $N$, its definition does not require any auxiliary structures like metrics and connections.}\adjustfootnotemark{-1}
    Let $\bold T[R]$ be $R$-th native extended tensor field associated
with the tensor bundle $N=T^{r_1\ldots\,r_Q}_{s_1\ldots\,s_Q}\!M$ and
let $\bold Y$ be some arbitrary extended tensor field of the type $(r_P,
s_P)$. Then we can apply $\vnabla_{\!\bold Y}[P]$ to $\bold T[R]$. By
means of the direct calculations using the explicit formula 
\mythetag{12.11} in local coordinates we find that
$$
\vnabla_{\!\bold Y}[P]\bold T[R]
=\cases \bold Y &\text{for \ }P=R,\\
0 &\text{for \ }P\neq R.\endcases
$$
Like covariant differentiations (see theorem~\mythetheorem{11.1} and
definition~\mythedefinition{11.4}), multivariate differentiations are
associated with some lifts. However, unlike covariant differentiations,
they lift not vectors, but tensors, though converting them into tangent
vectors of the bundle $N$. \pagebreak In the case of the canonical multivariate differentiation $\vnabla[P]$ for each point $q\in N$ we have some 
$\Bbb R$-linear map
$$
\hskip -2em
f[P]\!:\,T^r_s(\pi(q),M)\to T_q(N),
\mytag{12.12}
$$ 
where $r=r_P$ and $s=s_P$. The map \mythetag{12.12} takes a tensor $\bold
Y\in T^r_s(\pi(q),M)$ to the following vector in the tangent space $T_q(N)$
of the manifold $N$ at the point $q$:
$$
\hskip -2em
f[P](\bold Y)=
\dsize\msum{n}\Sb h_1,\,\ldots,\,h_r\\k_1,\,\ldots,\,k_s\endSb
Y^{h_1\ldots\,h_r}_{k_1\ldots\,k_s}\
\bold V^{k_1\ldots\,k_s}_{h_1\ldots\,h_r}[P].
\mytag{12.13}
$$
Here again $r=r_P$, $s=s_P$, and the vectors $\bold V^{k_1\ldots\,
k_s}_{h_1\ldots\,h_r}[P]$ are given by the second formula \mythetag{3.7}.
In a coordinate-free form the formula \mythetag{12.13} can be interpreted
as follows: the vector $f[P](\bold Y)$ in \mythetag{12.13} is the tangent
vector of the parametric curve \mythetag{12.5} at its initial point 
$q=q(0)$.\par
     Let's consider the image of the $\Bbb R$-linear map \mythetag{12.12}.
We denote it $V_q[P](N)$. Then from \mythetag{12.13} one easily derives 
that $V_q[P](N)$ is a subspace within the vertical subspace $V_q(N)$ of 
the tangent space $T_q(N)$. Moreover, we have
$$
\hskip -2em
V_q(N)=V_q[1](N)\oplus\ldots\oplus V_q[Q](N).
\mytag{12.14}
$$
The formula \mythetag{12.14} is a well-known fact, it follows from
\mythetag{3.1}. Due to the inclusion
$$
\Img f[P]=V_q[P](N)\subset V_q(N)
$$
the multivariate differentiation \mythetag{12.11} is a vertical
differentiation.\par
\head
13. Horizontal covariant differentiations\\
and extended connections.
\endhead
    Let $\nabla$ be some horizontal covariant differentiation of 
the algebra of extended tensor fields $\bold T(M)$ and let $f$ be 
the horizontal lift of vectors associated with it (see 
definition~\mythedefinition{11.4}). The horizontality of $f$ means 
that the image of the linear map \mythetag{11.2} is some $n$-dimensional
subspace $H_q(N)$ within the tangent space $T_q(N)$. It is called a
{\it horizontal subspace}. Due to $\pi_*\compos f=\idop$ the mappings
$$
\xalignat 2
&\hskip -2em
f\!:T_{\pi(q)}(M)\to H_q(N),
&&\pi_*\!:H_q(N)\to T_{\pi(q)}(M)
\mytag{13.1}
\endxalignat
$$
are inverse to each other. Due to the same equality $\pi_*\compos f=\idop$
the sum of the vertical and horizontal subspaces is a direct sum:
$$
\hskip -2em
H_q(N)\oplus V_q(N)=T_q(N).
\mytag{13.2}
$$
\mytheorem{13.1} Defining a horizontal lift of vectors from $M$ to $N$
is equivalent to fixing some direct complement $H_q(N)$ of the vertical
subspace $V_q(N)$ within the tangent space $T_q(N)$ at each point 
$q\in N$.
\endproclaim
\demo{Proof} Suppose that some horizontal lift of vectors $f$ is given.
\pagebreak 
Then the subspace $H_q(N)$ at the point $q$ is determined as the image
of the mapping \mythetag{11.2}, while the relationship \mythetag{13.2}
is derived from $\pi_*\compos f=\idop$ and from \mythetag{11.3}.\par
     Conversely, assume that at each point $q\in N$ we have a subspace
$H_q(N)$ complementary to $V_q(N)$. Then at each point $q\in N$ the
relationship \mythetag{13.2} is fulfilled. The kernel of the mapping
$\pi_*\!:T_q(N)\to T_{\pi(q)}(M)$ coincides with $V_q(N)$, therefore the
restriction of $\pi_*$ to the horizontal subspace $H_q(N)$ is a bijection.
The lift of vectors $f$ from $M$ to $N$ then can be defined as the inverse
mapping for $\pi_*\!:H_q(N)\to T_{\pi(q)}(M)$. If $f$ is defined in this 
way, then the mappings \mythetag{13.1} appear to be inverse to each other
and we get the equality $\pi_*\compos f=\idop$. According to the
definition~\mythedefinition{11.3}, it means that $f$ is a horizontal mapping. The theorem is completely proved.
\qed\enddemo
    Let's study a horizontal lift of vectors $f$ in a coordinate form. 
Upon choosing some local chart in $M$ we can consider the coordinate 
vector fields $\bold E_1,\,\ldots,\,\bold E_n$ in this chart (see \mythetag{1.3}). Applying the lift $f$ to them, we get
$$
\hskip -2em
f(\bold E_j)=\bold U_j-\sum^Q_{P=1}
\dsize\msum{n}\Sb i_1,\,\ldots,\,i_r\\j_1,\,\ldots,\,j_s\endSb
\Gamma^{\ i_1\ldots\,i_r}_{j\,j_1\ldots\,j_s}[P]\ \bold V^{j_1\ldots
\,j_s}_{i_1\ldots\,i_r}[P].
\mytag{13.3}
$$
Here $r=r_P$ and $s=s_P$, while $\bold U_j$ and $\bold V^{j_1\ldots
\,j_s}_{i_1\ldots\,i_r}[P]$ are determined by \mythetag{3.7} This
formula for $f(\bold E_j)$ follows from $\pi_*\compos f=\idop$ due 
to the equalities 
$$
\xalignat 2
&\hskip -2em
\pi_*\!\left(\bold U_j\right)=\bold E_j,
&&\pi_*\!\left(\bold V^{j_1\ldots\,j_s}_{i_1\ldots\,i_r}[P]\right)=0.
\mytag{13.4}
\endxalignat
$$
The quantities $\Gamma^{\ i_1\ldots\,i_r}_{j\,j_1\ldots\,j_s}[P]$ in 
\mythetag{13.3} are called the {\it components} of a horizontal lift 
of vectors in a local chart $U$. If the lift $f$ is induced by some
horizontal covariant differentiation $\nabla$, then for its components
in\mythetag{13.4} we have
$$
\hskip -2em
\Gamma^{\ i_1\ldots\,i_r}_{j\,j_1\ldots\,j_s}[P]=
-Z^{\ i_1\ldots\,i_r}_{j\,j_1\ldots\,j_s}[P].
\mytag{13.5}
$$
The quantities $Z^{\ i_1\ldots\,i_r}_{j\,j_1\ldots\,j_s}[P]$ in  \mythetag{13.5} are the same as in \mythetag{9.3}, \mythetag{9.5}, \mythetag{9.6}, and in \mythetag{9.7}. As for the quantities $Z^i_j$ 
in \mythetag{9.4}, in the case of a horizontal covariant differentiation
they are given by the Kronecker's delta-symbol: $Z^i_j=\delta^i_j$.
Substituting $Z^i_j=\tilde Z^i_j=\delta^i_j$ into \mythetag{9.6} and 
taking into account \mythetag{13.5}, we derive
$$
\pagebreak
\align
&\cases
\gathered
\Gamma^{\ i_1\ldots\,i_r}_{j\,j_1\ldots\,j_s}[P]=\sum^n_{k=1}
\dsize\msum{n}\Sb h_1,\,\ldots,\,h_r\\k_1,\,\ldots,\,k_s\endSb 
S^{i_1}_{h_1}\ldots\,S^{i_r}_{h_r}
\ T^{k_1}_{j_1}\ldots\,T^{k_s}_{j_s}\,T^k_j\ 
\tilde\Gamma^{\,h_1\ldots\,h_r}_{k\,k_1\ldots\,k_s}[P]\ +\\
+\sum^r_{m=1}\sum^n_{v_m=1}
\theta^{i_m}_{\!j\,v_m}\,T^{\,i_1\ldots\,v_m\ldots\,i_r}_{j_1
\ldots\,\ldots\,\ldots\,j_s}[P]
-\sum^s_{m=1}\sum^n_{w_m=1}
\theta^{w_m}_{\!j\,j_m}\,T^{\,i_1\ldots\,\ldots\,\ldots\,
i_r}_{j_1\ldots\,w_m\ldots\,j_s}[P],
\endgathered\endcases\ \qquad
\vspace{-45pt}
\mytag{13.6}\\
\vspace{35pt}
&\cases
\gathered
\tilde\Gamma^{\ i_1\ldots\,i_r}_{j\,j_1\ldots\,j_s}[P]=\sum^n_{k=1}
\dsize\msum{n}\Sb h_1,\,\ldots,\,h_r\\k_1,\,\ldots,\,k_s\endSb 
T^{i_1}_{h_1}\ldots\,T^{i_r}_{h_r}
\ S^{k_1}_{j_1}\ldots\,S^{k_s}_{j_s}\,S^k_j\ 
\Gamma^{\,h_1\ldots\,h_r}_{k\,k_1\ldots\,k_s}[P]\ +\\
+\sum^r_{m=1}\sum^n_{v_m=1}
\tilde\theta^{i_m}_{\!j\,v_m}\,\tilde T^{\,i_1\ldots\,v_m\ldots\,i_r}_{j_1
\ldots\,\ldots\,\ldots\,j_s}[P]
-\sum^s_{m=1}\sum^n_{w_m=1}
\tilde\theta^{w_m}_{\!j\,j_m}\,\tilde T^{\,i_1\ldots\,\ldots\,\ldots\,
i_r}_{j_1\ldots\,w_m\ldots\,j_s}[P],
\endgathered
\endcases\ \qquad
\vspace{-40pt}
\mytag{13.7}
\endalign
$$
Here $r=r_P$, $s=s_P$. The $\theta$-parameters are taken from
\mythetag{3.15}, \mythetag{3.10}, and \mythetag{3.12}.
The formulas \mythetag{13.6} and \mythetag{13.7} express the 
transformation rules for the components of a horizontal lift 
of vectors in \mythetag{13.3}.\par
     Another geometric structure associated with a horizontal covariant
differentiation $\nabla$ reveals when we apply $D=\nabla_{\bold Y}$ to
the module $T^1_0(M)$. This produces the mapping \mythetag{7.23}. In a
local chart it is described by the formula \mythetag{7.13}. In the
present case we can take $\bold Y=\hat\bold E_j$ and write this formula
as
$$
\hskip -2em
\nabla_{\hat\bold E_j}\hat\bold E_i=\sum^n_{k=1}\Gamma^k_{j\,i}\,\bold E_k.
\mytag{13.8}
$$
Here $\hat\bold E_i$ and $\hat\bold E_j$ are determined by the 
formula \mythetag{7.13}. Relying on the lemma~\mythelemma{7.4} 
and on the localization theorem~\mythetheorem{8.1}, we can write 
\mythetag{13.8} as follows:
$$
\hskip -2em
\nabla_{\bold E_j}\bold E_i=\sum^n_{k=1}\Gamma^k_{j\,i}\,\bold E_k.
\mytag{13.9}
$$
The coefficients $\Gamma^k_{j\,i}$ are the same as in \mythetag{9.4}.
Since $Z^i_j=\delta^i_j$ for a horizontal covariant differentiation,
the transformation formulas \mythetag{9.8} and \mythetag{9.9} for the
coefficients $\Gamma^k_{j\,i}$ in \mythetag{13.8} and \mythetag{13.9}
now reduce to the following ones:
$$
\align
&\hskip -2em
\Gamma^k_{j\,i}=\dsize\sum^n_{b=1}\sum^n_{a=1}\sum^n_{c=1}
S^k_a\,T^b_i\,T^c_j\ \tilde\Gamma^a_{c\,b}+
\sum^n_{a=1}\theta^k_{j\,i},
\mytag{13.10}\\
\vspace{2ex}
&\hskip -2em
\tilde\Gamma^k_{j\,i}=\dsize\sum^n_{b=1}\sum^n_{a=1}\sum^n_{c=1}
T^k_a\,S^b_i\,S^c_j\ \Gamma^a_{c\,b}+
\sum^n_{a=1}\tilde\theta^k_{j\,i}.
\mytag{13.11}
\endalign
$$
\mydefinition{13.1} Let $M$ be a smooth manifold and let
$N=T^{r_1\ldots\,r_Q}_{s_1\ldots\,s_Q}\!M$ be a composite tensor
bundle over $M$. An extended affine connection $\Gamma$ is a geometric
object in each local chart of $M$ represented by its components
$\Gamma^k_{j\,i}$ and such that its components are smooth functions
of the variables \mythetag{3.3} transformed according to the formulas
\mythetag{13.10} and \mythetag{13.11} under a change of a local 
chart.
\enddefinition
\mytheorem{13.2} On any smooth paracompact manifold $M$ equipped with
a composite tensor bundle $N=T^{r_1\ldots\,r_Q}_{s_1\ldots\,s_Q}\!M$ 
there is at least one extended affine connection.
\endproclaim
     We shall not prove this theorem here. Its proof for the spacial
case, where $N=TM$, is given in Chapter~\uppercase
\expandafter{\romannumeral 3} of the thesis \mycite{4}. This proof 
can be easily transformed for the present more general case. Note also
that any traditional affine connection fits the above  definition~\mythedefinition{13.1} being a special case for this more 
general concept of an extended connection.
\mydefinition{13.2} A horizontal covariant differentiation $\nabla$ 
of the algebra of extended tensor fields $\bold T(M)$ associated with 
some composite tensor bundle $T^{r_1\ldots\,r_Q}_{s_1\ldots\,s_Q}\!M$
is called a {\it spatial covariant differentiation} or a {\it spatial
gradient} if
$$
\hskip -2em
\nabla\bold T[P]=0\text{ \ \ for all \ \ }P=1,\,\ldots,\,Q,
\mytag{13.12}
$$
i\.\,e\. if the operator $\nabla$ annuls all native extended tensor fields $\bold T[1],\,\ldots,\,\bold T[Q]$.
\enddefinition
    Let's study the equality \mythetag{13.12} specifying spatial 
covariant differentiations. For this purpose we use the formula
\mythetag{9.5} substituting $Z^i_j=\delta^i_j$ and \mythetag{13.5}
into it:
$$
\gathered
\nabla_{\!\!j}X^{i_1\ldots\,i_\alpha}_{j_1\ldots\,j_\beta}
=\frac{\partial X^{i_1\ldots\,i_\alpha}_{j_1\ldots\,j_\beta}}
{\partial x^j}-\sum^Q_{R=1}\dsize\msum{n}\Sb h_1,\,\ldots,\,h_r\\k_1,\,\ldots,\,k_s\endSb \Gamma^{\,h_1\ldots\,h_r}_{j\,k_1\ldots\,k_s}[R]
\ \frac{\partial X^{i_1\ldots\,i_\alpha}_{j_1\ldots\,j_\beta}}
{\partial T^{h_1\ldots\,h_r}_{k_1\ldots\,k_s}[R]
\vphantom{\vrule height 10pt depth 0pt}}\ +\\
+\sum^\alpha_{m=1}\sum^n_{v_m=1}\Gamma^{i_m}_{j\,v_m}\ X^{\,i_1\ldots\,v_m
\ldots\,i_\alpha}_{j_1\ldots\,\ldots\,\ldots\,j_\beta}-\sum^\beta_{m=1}
\sum^n_{w_m=1}\Gamma^{w_m}_{j\,j_m}\ X^{\,i_1\ldots\,\ldots\,\ldots\,
i_\alpha}_{j_1\ldots\,w_m\ldots\,j_\beta},
\endgathered\quad
\mytag{13.13}
$$
According to the formula \mythetag{13.12}, we should substitute
$\alpha=r=r_P$, $\beta=s=s_P$, and $X^{i_1\ldots\,i_\alpha}_{j_1\ldots
\,j_\beta}=T^{i_1\ldots\,i_r}_{j_1\ldots\,j_s}[P]$ into the formula
\mythetag{13.13}. Recall that the quantities $T^{i_1\ldots\,i_r}_{j_1
\ldots\,j_s}[P]$ in \mythetag{9.5} and \mythetag{13.13} are treated as
independent variables. Therefore, we get
$$
\gathered
\nabla_{\!\!j}T^{i_1\ldots\,i_r}_{j_1\ldots\,j_s}[P]=
-\Gamma^{\,i_1\ldots\,i_r}_{j\,j_1\ldots\,j_s}[P]\ +\\
+\ \sum^r_{m=1}\sum^n_{v_m=1}\Gamma^{i_m}_{j\,v_m}\ T^{\,i_1\ldots\,v_m
\ldots\,i_r}_{j_1\ldots\,\ldots\,\ldots\,j_s}[P]-\sum^s_{m=1}
\sum^n_{w_m=1}\Gamma^{w_m}_{j\,j_m}\ T^{\,i_1\ldots\,\ldots\,\ldots\,
i_r}_{j_1\ldots\,w_m\ldots\,j_s}[P]=0.
\endgathered
$$
This formula can be rewritten in the following form:
$$
\hskip -2em
\gathered
\Gamma^{\,i_1\ldots\,i_r}_{j\,j_1\ldots\,j_s}[P]=\sum^r_{m=1}
\sum^n_{v_m=1}\Gamma^{i_m}_{j\,v_m}\ T^{\,i_1\ldots\,v_m
\ldots\,i_r}_{j_1\ldots\,\ldots\,\ldots\,j_s}[P]\ -\\
-\sum^s_{m=1}\sum^n_{w_m=1}\Gamma^{w_m}_{j\,j_m}\
T^{\,i_1\ldots\,\ldots\,\ldots\,i_r}_{j_1\ldots\,w_m\ldots\,j_s}[P].
\endgathered
\mytag{13.14}
$$
Substituting \mythetag{13.14} back into \mythetag{13.13}, we derive
$$
\gathered
\nabla_{\!\!j}X^{i_1\ldots\,i_\alpha}_{j_1\ldots\,j_\beta}
=\frac{\partial X^{i_1\ldots\,i_\alpha}_{j_1\ldots\,j_\beta}}
{\partial x^j}\,+\sum^\alpha_{m=1}\sum^n_{v_m=1}\Gamma^{i_m}_{\!j
\,v_m}\ X^{\,i_1\ldots\,v_m\ldots\,i_\alpha}_{j_1\ldots\,\ldots
\,\ldots\,j_\beta}\ -\\
-\sum^Q_{R=1}
\dsize\msum{n}\Sb h_1,\,\ldots,\,h_r\\k_1,\,\ldots,\,k_s\endSb
\sum^r_{m=1}\sum^n_{v_m=1}\Gamma^{h_m}_{\!j\,v_m}\ T^{\,h_1\ldots\,v_m
\ldots\,h_r}_{\ k_1\ldots\,\ldots\,\ldots\,k_s}[R]\ \frac{\partial X^{i_1
\ldots\,i_\alpha}_{j_1\ldots\,j_\beta}}{\partial T^{h_1\ldots\,h_r}_{k_1
\ldots\,k_s}[R]\vphantom{\vrule height 10pt depth 0pt}}\ -\\
-\sum^\beta_{m=1}
\sum^n_{w_m=1}\Gamma^{w_m}_{\!j\,j_m}\ X^{\,i_1\ldots\,\ldots\,\ldots\,
i_\alpha}_{j_1\ldots\,w_m\ldots\,j_\beta}\ +\\
+\sum^Q_{R=1}
\dsize\msum{n}\Sb h_1,\,\ldots,\,h_r\\k_1,\,\ldots,\,k_s\endSb
\sum^s_{m=1}\sum^n_{w_m=1}\Gamma^{w_m}_{\!j\,k_m}\
T^{\,h_1\ldots\,\ldots\,\ldots\,h_r}_{k_1\ldots\,w_m\ldots\,k_s}[R]
\ \frac{\partial X^{i_1\ldots\,i_\alpha}_{j_1\ldots\,j_\beta}}
{\partial T^{h_1\ldots\,h_r}_{k_1\ldots\,k_s}[R]
\vphantom{\vrule height 10pt depth 0pt}}.
\endgathered\quad
\mytag{13.15}
$$\par
The theorem~\mythetheorem{7.2} applied to a horizontal covariant
differentiation says that any such differentiation is defined by
two independent geometric structures: 
\roster
\item a horizontal lift of vectors from $M$ to $N$;
\item an extended connection.
\endroster
The formula \mythetag{13.14} relates these two structures. It 
expresses the components of the horizontal lift in \mythetag{13.3} 
through the components of an extended connection $\Gamma$ in \mythetag{13.9}. \pagebreak This result is formulated as the 
following theorem.
\mytheorem{13.3} Defining a spacial covariant differentiation in the
algebra of extended tensor fields $\bold T(M)$ is equivalent to defining
an extended connection $\Gamma$.
\endproclaim
\head
14. The structural theorem for differentiations.
\endhead
\mytheorem{14.1} Let $M$ be a smooth manifold and let $N=T^{r_1\ldots
\,r_Q}_{s_1\ldots\,s_Q}\!M$ be a composite tensor bundle over $M$. 
If\/ $M$ is equipped with some extended affine connection $\Gamma$, 
then each differentiation $D$ of the algebra of extended tensor fields 
$\bold T(M)$ in this manifold $M$ is uniquely expanded into a sum
$$
\hskip -2em
D=\nabla_{\bold X}+\sum^Q_{P=1}\vnabla_{\bold Y_{\!P}}[P]+S,
\mytag{14.1}
$$
where $\nabla_{\bold X}$ is the spacial covariant differentiation along
some extended vector field $\bold X$, $\vnabla_{\bold Y_{\!P}}[P]$ is 
the $P$-th canonical vertical multivariate differentiation along some 
extended tensor field $\bold Y_{\!P}$ of the type $(r_P,s_P)$, and 
$S$ is a degenerate differentiation given by some extended tensor
field $\bold S$ of the type $(1,1)$.
\endproclaim
\demo{Proof} Let $D\in\goth D(M)$. Then its restriction to $T^0_0(M)$
is given by some vector field $\bold Z$ in $N$. The extended affine
connection $\Gamma$ in $M$ determines some horizontal lift of vectors $f$
from $M$ to $N$. Its components in a local chart are given by the formula
\mythetag{13.14}. According to the theorem~\mythetheorem{13.1}, this lift
of vectors determines the expansion of the tangent space $T_q(N)$ into
a direct sum \mythetag{13.2} at each point $q\in N$. The vertical subspace
$V_q(N)$ in \mythetag{13.2} has its own expansion \mythetag{12.14} into a
direct sum. Combining \mythetag{13.2} and \mythetag{12.14}, we obtain
$$
\hskip -2em
T_q(N)=H_q(N)\oplus V_q[1](N)\oplus\ldots\oplus V_q[Q](N).
\mytag{14.2}
$$
Then the vector field $\bold Z$ is expanded into a sum of vector fields
$$
\hskip -2em
\bold Z=\bold H+\bold V_{\!1}+\ldots+\bold V_{\!Q}
\mytag{14.3}
$$
uniquely determined by the expansion \mythetag{14.2}. Due to the maps 
\mythetag{12.12} each vector field $\bold V_{\!P}$ in \mythetag{14.3}
is uniquely associated with some extended tensor field $\bold Y_{\!P}$
of the type $(r_P,s_P)$. Similarly, the vector field $\bold H$ is 
uniquely associated with the the extended vector field $\bold X$ such 
that $\bold H_q=f(\bold X_q)$. Then we can consider the sum
$$
\hskip -2em
\tilde D=\nabla_{\bold X}+\sum^Q_{P=1}\vnabla_{\bold Y_{\!P}}[P].
\mytag{14.4}
$$
The sum \mythetag{14.4} is a differentiation of $\bold T(M)$ such that
its restriction to $T^0_0(M)$ is given by the vector \mythetag{14.3}.
Hence, $D-\tilde D$ is a differentiation of $\bold T(M)$ with identically
zero restriction to $T^0_0(M)$. This means that $D-\tilde D$ is a
degenerate differentiation (see definition~\mythedefinition{8.1}). 
Applying the theorem~\mythetheorem{8.3}, we find that $S=D-\tilde D$ is
given by some extended tensor field $\bold S$ of the type $(1,1)$. Thus,
the expansion \mythetag{14.1} and the theorem~\mythetheorem{14.1} 
in whole are proved.\qed\enddemo
    The theorem~\mythetheorem{14.1} is the structural theorem for
differentiations in the algebra of extended tensor fields $\bold T(M)$.
It approves our previous efforts in studying the three basic types of
differentiations which are used in the formula \mythetag{14.1}.
\head
\S\,15. Commutation relationships and curvature tensors.
\endhead
    Let's remember that the set of all differentiations $\goth D(M)$
is an infinite-dimensional Lie algebra (see formula \mythetag{6.1}).
Using the above structural theorem~\mythetheorem{14.1}, one can give
a more detailed description of this Lie algebra. Let's begin with 
degenerate differentiations. Assume that $S_1$ and $S_2$
are two degenerate differentiations given by two extended tensor fields
$\bold S_1$ and $\bold S_2$ of the type $(1,1)$. Then
$$
\hskip -2em
[S_1,\,S_2]=S_3\text{, \ wehere \ }\bold S_3=C(\bold S_1\otimes
\bold S_2-\bold S_2\otimes\bold S_1).
\mytag{15.1}
$$
The formula \mythetag{15.1} means that the commutator of two degenerate
differentiations is a degenerate differentiation given by the pointwise
commutator of the corresponding extended tensor fields $\bold S_1$ and 
$\bold S_2$.\par
     Assume that $M$ is equipped with an extended affine connection $\Gamma$. Then we can consider the commutators of some degenerate differentiation $\bold S$ with the spacial covariant differentiation $\nabla_{\!\bold X}$ and with the $P$-th canonical vertical 
multivariate differentiation $\vnabla_{\bold Y}[P]$. These commutators
are given by the formula
$$
\alignat 2
&\hskip -2em
[\nabla_{\!\bold X},\,S]=S_1,
&\quad &\text{where \ }\bold S_1=\nabla_{\!\bold X}\bold S;\\
\vspace{-1ex}
\mytag{15.2}\\
\vspace{-1ex}
&\hskip -2em
[\vnabla_{\bold Y}[P],\,\bold S]=S_2,
&\quad &\text{where \ }\bold S_2=\vnabla_{\bold Y}[P]\bold S.
\endalignat
$$
The formulas \mythetag{15.2} mean that both commutators are again
degenerate differentiations. They are given by the extended tensor fields
$\nabla_{\!\bold X}\bold S$ and $\vnabla_{\bold Y}[P]\bold S$ respectively.
\par
    The commutator of two canonical vertical multivariate differentiations 
$\vnabla_{\!\bold X}[P]$ and $\vnabla_{\bold Y}[R]$ is composed by other two
such differentiations. Indeed, we have
$$
\hskip -2em
[\vnabla_{\!\bold X}[P],\,\vnabla_{\bold Y}[R]]
=\vnabla_{\!\bold U}[R]-\vnabla_{\bold V}[P].
\mytag{15.3}
$$
where $\bold U$ and $\bold V$ are determined as follows:
$$
\xalignat 2
&\hskip -2em
\bold U=\vnabla_{\!\bold X}[P]\bold Y,
&&\bold V=\vnabla_{\bold Y}[R]\bold X.
\mytag{15.4}
\endxalignat
$$
Similarly, for the commutator of the spatial covariant
differentiation $\nabla_{\!\bold X}$ with the canonical vertical multivariate differentiation $\vnabla_{\bold Y}[P]$ we get
$$
\hskip -2em
[\nabla_{\!\bold X},\,\vnabla_{\bold Y}[P]]=
\vnabla_{\!\bold U}[P]+\sum^Q_{R=1}
\vnabla_{\!\bold U[R]}[R]-\nabla_{\bold V}+S,
\mytag{15.5}
$$
where $\bold U$, $\bold U[R]$, and $\bold V$ are determined as follows:
$$
\xalignat 3
&\hskip -2em
\bold U=\nabla_{\!\bold X}\bold Y,
&&\bold U[R]=-S\,\bold T[R],
&&\bold V=\vnabla_{\bold Y}[P]\bold X.\qquad
\mytag{15.6}
\endxalignat
$$
As for $S$ in \mythetag{15.5} and \mythetag{15.6}, it is a degenerate
differentiation determined by some definite extended tensor field $\bold S$
of the type $(1,1)$ depending on $\bold X$ and on $\bold Y$:
$$
\hskip -2em
\bold S=\bold D[P](\bold X,\bold Y)=C(\bold D[P]\otimes\bold X\otimes
\bold Y).
\mytag{15.7}
$$ 
Similarly, $\bold U[R]$ in \mythetag{15.6} is some definite \pagebreak
extended tensor field of the type $(r_P,s_P)$ depending on $\bold X$, 
on $\bold Y$, and on the indices $P$ and $R$: 
$$
\hskip -2em
\bold U[R]
=\boldsymbol\Theta[P,R](\bold X,\bold Y)=C(\boldsymbol\Theta[P,R]\otimes
\bold X\otimes\bold Y).
\mytag{15.8}
$$
The basic object in the series of notations \mythetag{15.6},
\mythetag{15.7}, \mythetag{15.8} is $\bold D[P]$. It is called the
{\it $P$-th dynamic curvature tensor}. This is an extended tensor 
field of the type $(s_P+1,r_P+2)$. Its components in a local chart
are given by the formula
$$
\hskip -2em
D^{\ k\,k_1\ldots\,k_s}_{ij\,h_1\,\ldots\,h_r}[P]=
-\frac{\partial\Gamma^k_{j\,i}}{\partial 
T^{h_1\,\ldots\,h_r}_{k_1\ldots\,k_s}[P]},
\mytag{15.9}
$$
where $r=r_P$ and $s=s_P$. Then in \mythetag{15.8} we have the extended
tensor field $\boldsymbol\Theta[P,R]$ of the type $(r_R+s_P,s_R+r_P+1)$.
Its components are expressed through the components of $\bold D[P]$ in
\mythetag{15.9} according to the formula
$$
\gathered
\Theta^{\ i_1\ldots\,i_\alpha k_1\ldots\,k_s}_{j_1\ldots\,j_\beta\,j\, h_1\,\ldots\,h_r}[P,R]=\sum^\beta_{m=1}\sum^n_{w_m=1}D^{w_m\,h_1
\,\ldots\,h_r}_{j_m\,j\,k_1\ldots\,k_s}[P]\ T^{\,i_1\ldots\,\ldots\,
\ldots\,i_\alpha}_{j_1\ldots\,w_m\ldots\,j_\beta}[R]\ -\\
-\sum^\alpha_{m=1}\sum^n_{v_m=1}D^{i_m\,h_1\,\ldots\,h_r}_{v_m\,j
\,k_1\ldots\,k_s}[P]\ T^{\,i_1\ldots\,v_m\ldots\,i_\alpha}_{j_1\ldots
\,\ldots\,\ldots\,j_\beta}[R],
\endgathered\quad
\mytag{15.10}
$$
where $r=r_P$, $s=s_P$, $\alpha=r_R$, and $\beta=s_R$. The components of
the tensor $\bold U[R]=\boldsymbol\Theta[P,R](\bold X,\bold Y)$
in \mythetag{15.8} are expressed through \mythetag{15.10} as follows:
$$
U^{\,i_1\ldots\,i_\alpha}_{j_1\ldots\,j_\beta}[R]=
\sum^n_{j=1}\dsize\msum{n}\Sb h_1,\,\ldots,\,h_r\\k_1,\,\ldots,\,k_s
\endSb\Theta^{\ i_1\ldots\,i_\alpha k_1\ldots\,k_s}_{j_1\ldots\,
j_\beta\,j\,h_1\,\ldots\,h_r}[P,R]\ X^j\ Y^{h_1\,\ldots\,h_r}_{k_1\ldots
\,k_s}.\quad
\mytag{15.11}
$$
Similarly, the components of the tensor $\bold S=\bold D[P](\bold X,
\bold Y)$ in \mythetag{15.7} are expressed through \mythetag{15.9}
according to the following formula:
$$
\hskip -2em
S^k_i=\sum^n_{j=1}\dsize\msum{n}\Sb h_1,\,\ldots,\,h_r\\k_1,\,\ldots,\,
k_s\endSb D^{\ k\,k_1\ldots\,k_s}_{ij\,h_1\,\ldots\,h_r}[P]\ X^j\ Y^{h_1
\,\ldots\,h_r}_{k_1\ldots\,k_s}.
\mytag{15.12}
$$
The formula \mythetag{15.10} is derived from the second formula \mythetag{15.6} due to \mythetag{15.7} and \mythetag{15.8}. The 
formulas \mythetag{15.11} and \mythetag{15.12} are rather 
obvious. They complete the series of equalities which are used 
in order to make certain the right hand side of the commutation relationship \mythetag{15.5}.\par
    In the last step now we consider the commutator of two spatial
covariant differentiations $\nabla_{\!\bold X}$ and $\nabla_{\bold Y}$.
The formula for this commutator is written as 
$$
\hskip -2em
[\nabla_{\!\bold X},\,\nabla_{\bold Y}]=
\nabla_{\!\bold U}+\sum^Q_{R=1}
\vnabla_{\!\bold U[P]}[P]+S,
\mytag{15.13}
$$
where $\bold U$ and $\bold U[R]$ are determined in the following way:
$$
\xalignat 2
&\hskip -2em
\bold U=\nabla_{\!\bold X}\bold Y-\nabla_{\bold Y}\bold X-\bold V,
&&\bold U[R]=-S\,\bold T[R].\qquad
\mytag{15.14}
\endxalignat
$$
Like in \mythetag{15.5}, by $S$ in \mythetag{15.13} and 
\mythetag{15.14} we denote a degenerate differentiation 
determined by some definite extended tensor field $\bold S$
of the type $(1,1)$ depending on both extended vector fields 
$\bold X$ and $\bold Y$:
$$
\hskip -2em
\bold S=\bold R(\bold X,\bold Y)=C(\bold R\otimes\bold X
\otimes\bold Y).
\mytag{15.15}
$$
Similarly, $\bold V$ in \mythetag{15.6} is some definite extended 
vector field depending on $\bold X$ and on $\bold Y$. It is expressed 
through the {\it torsion tensor\/} $\bold T$:
$$
\hskip -2em
\bold V=\bold T(\bold X,\bold Y)=C(\bold T
\otimes\bold X\otimes\bold Y).
\mytag{15.16}
$$
The components of the torsion tensor in a local chart are given by the
formula
$$
\hskip -2em
T^k_{ij}=\Gamma^k_{ij}-\Gamma^k_{j\,i}.
\mytag{15.17}
$$
This formula \mythetag{15.17} coincides with the standard formula for
torsion (see \mycite{35}). The only difference here is that $\Gamma$ is assumed to be an extended connection, therefore $\bold T$ is an extended
tensor field of the type $(1,2)$.\par
    For the parameter $\bold U[R]$ in \mythetag{15.14} we write the 
formula analogous to \mythetag{15.8} since this is also some definite
extended tensor field depending on $\bold X$ and $\bold Y$:
$$
\hskip -2em
\bold U[R]=\boldsymbol\Omega[R](\bold X,\bold Y)=C(\boldsymbol\Omega[R]
\otimes\bold X\otimes\bold Y).
\mytag{15.18}
$$
The basic object in the series of notations \mythetag{15.14},
\mythetag{15.15}, \mythetag{15.16}, and \mythetag{15.18} is the
{\it curvature tensor\/} $\bold R$. In contrast to $\bold D[P]$
in \mythetag{15.9}, we call it the {\it static curvature tensor}.
For the components of the static curvature tensor we have the formula
$$
\gathered
R^k_{h\,ij}=\frac{\partial\Gamma^k_{jh}}{\partial x^i}
-\frac{\partial\Gamma^k_{\vphantom{j}ih}}{\partial x^j}
+\sum^n_{a=1}\Gamma^a_{jh}\,\Gamma^k_{i\,a}-\sum^n_{a=1}
\Gamma^a_{ih}\,\Gamma^k_{j\,a}-\\
-\sum^Q_{P=1}
\dsize\msum{n}\Sb h_1,\,\ldots,\,h_r\\k_1,\,\ldots,\,k_s\endSb
\sum^r_{m=1}\sum^n_{v_m=1}\Gamma^{h_m}_{\!i\,v_m}\ T^{\,h_1\ldots\,v_m
\ldots\,h_r}_{\ k_1\ldots\,\ldots\,\ldots\,k_s}[P]\ \frac{\partial 
\Gamma^k_{jh}}{\partial T^{h_1\ldots\,h_r}_{k_1
\ldots\,k_s}[P]\vphantom{\vrule height 10pt depth 0pt}}\ +\\
+\sum^Q_{P=1}
\dsize\msum{n}\Sb h_1,\,\ldots,\,h_r\\k_1,\,\ldots,\,k_s\endSb
\sum^r_{m=1}\sum^n_{v_m=1}\Gamma^{h_m}_{\!j\,v_m}\ T^{\,h_1\ldots\,v_m
\ldots\,h_r}_{\ k_1\ldots\,\ldots\,\ldots\,k_s}[P]\ \frac{\partial 
\Gamma^k_{ih}}{\partial T^{h_1\ldots\,h_r}_{k_1
\ldots\,k_s}[P]\vphantom{\vrule height 10pt depth 0pt}}\ +\\
+\sum^Q_{P=1}
\dsize\msum{n}\Sb h_1,\,\ldots,\,h_r\\k_1,\,\ldots,\,k_s\endSb
\sum^s_{m=1}\sum^n_{w_m=1}\Gamma^{w_m}_{\!i\,k_m}\
T^{\,h_1\ldots\,\ldots\,\ldots\,h_r}_{k_1\ldots\,w_m\ldots\,k_s}[P]
\ \frac{\partial\Gamma^k_{jh}}
{\partial T^{h_1\ldots\,h_r}_{k_1\ldots\,k_s}[P]
\vphantom{\vrule height 10pt depth 0pt}}\ -\\
-\sum^Q_{P=1}
\dsize\msum{n}\Sb h_1,\,\ldots,\,h_r\\k_1,\,\ldots,\,k_s\endSb
\sum^s_{m=1}\sum^n_{w_m=1}\Gamma^{w_m}_{\!j\,k_m}\
T^{\,h_1\ldots\,\ldots\,\ldots\,h_r}_{k_1\ldots\,w_m\ldots\,k_s}[P]
\ \frac{\partial\Gamma^k_{ih}}
{\partial T^{h_1\ldots\,h_r}_{k_1\ldots\,k_s}[P]
\vphantom{\vrule height 10pt depth 0pt}}.
\endgathered\qquad
\mytag{15.19}
$$
In the case of non-extended connection $\Gamma$ the formula 
\mythetag{15.19} reduces to the standard formula for the curvature
tensor (see \mycite{35}).\par
    Returning back to the equality \mythetag{15.18}, we need to 
write the formula for the components of the tensor $\boldsymbol
\Omega[R]$. Here is this formula:
$$
\gathered
\Omega^{\ i_1\ldots\,i_\alpha}_{\,ij\,j_1\ldots\,j_\beta }[R]=\sum^\beta_{m=1}\sum^n_{w_m=1}R^{w_m}_{j_mij}\ T^{\,i_1
\ldots\,\ldots\,\ldots\,i_\alpha}_{j_1\ldots\,w_m\ldots
\,j_\beta}[R]\ -\\
-\sum^\alpha_{m=1}\sum^n_{v_m=1}R^{i_m}_{v_m\,ij
}\ T^{\,i_1\ldots\,v_m\ldots\,i_\alpha}_{j_1\ldots
\,\ldots\,\ldots\,j_\beta}[R].
\endgathered\quad
\mytag{15.20}
$$
It is similar to \mythetag{15.10}. The formula \mythetag{15.20} is 
derived from \mythetag{15.18}, \mythetag{15.15}, and from the second 
formula \mythetag{15.14}. The analogs of the formulas \mythetag{15.11}
and \mythetag{15.12} in this case are written as follows:
$$
\gather
\hskip -2em
U^{\,i_1\ldots\,i_\alpha}_{j_1\ldots\,j_\beta}[R]=\sum^n_{i=1}
\sum^n_{j=1}\dsize\Omega^{\ i_1\ldots\,i_\alpha}_{\,ij\,j_1
\ldots\,j_\beta }[R]\ X^i\ Y^j,\quad
\mytag{15.21}\\
\hskip -2em
V^k=\sum^n_{i=1}\sum^n_{j=1}T^k_{ij}\ X^i\ Y^j,
\mytag{15.22}\\
\hskip -2em
S^k_h=\sum^n_{i=1}\sum^n_{j=1}R^k_{hij}\ X^i\ Y^j.
\mytag{15.23}
\endgather
$$
The formulas \mythetag{15.21}, \mythetag{15.22}, \mythetag{15.23}
complete the series of equalities which are written in order to make 
certain the right hand side of the commutation relationship \mythetag{15.13}. As for the commutation relationships themselves,
they can be derived by direct calculations on the base of the formulas
\mythetag{12.11} and \mythetag{13.15}.\par
\head
16. Coordinate representation of commutation relationships.
\endhead
     The first commutation relationship \mythetag{15.1} is trivial. In
coordinate form, i\.\,e\. in a local chart, it means that the matrix of
the tensor $\bold S_3$ is the matrix commutator produced from the matrices of $\bold S_1$ and $\bold S_2$.\par
     The next two commutator relationships \mythetag{15.2} are also
rather simple. They mean that the components of $\bold S_1$ and $\bold S_2$
are derived from the components of $\bold S$ by means of the formulas
\mythetag{12.11} and \mythetag{13.15}.\par
     The fourth commutation relationship \mythetag{15.3} is not so simple,
but in a coordinate form it reduces to the following one:
$$
\hskip -2em
[\vnabla^{j_1\ldots\,j_s}_{i_1\ldots\,i_r}[P],\ 
\vnabla^{k_1\ldots\,k_\beta}_{h_1\ldots\,h_\alpha}[R]]=0.
\mytag{16.1}
$$
The relationship \mythetag{16.1} is easily derived from \mythetag{12.11}.
\par
     Now let's proceed with the fifth commutation relationship
\mythetag{15.5}. In a local chart we should consider the commutator of
$\nabla_i$ and $\vnabla^{j_1\ldots\,j_s}_{i_1\ldots\,i_r}[P]$. From \mythetag{15.5} we derive
$$
\hskip -2em
\gathered
[\nabla_i,\,\vnabla^{j_1\ldots\,j_s}_{i_1\ldots\,i_r}[P]]X^k=
\sum^n_{h=1}D^{\,k\,j_1\ldots\,j_s}_{h\,i\,i_1\ldots\,i_r}[P]
\ X^h\ +\\
+\sum^Q_{R=1}\msum{n}\Sb h_1,\,\ldots,\,h_\alpha\\k_1,\,\ldots,
\,k_\beta\endSb
\Theta^{\,h_1\ldots\,h_\alpha\,j_1\ldots\,j_s}_{k_1\ldots
\,k_\beta\,i\,i_1\ldots\,i_r}[P,R]\ \vnabla^{k_1\ldots\,
k_\beta}_{h_1\ldots\,h_\alpha}[R]\ X^k.
\endgathered
\mytag{16.2}
$$
Here $X^1,\,\ldots,\,X^n$ are the components of some extended vector 
field $\bold X$. When applied to an extended scalar field $\varphi$ 
the same commutator is written as follows:
$$
[\nabla_i,\,\vnabla^{j_1\ldots\,j_s}_{i_1\ldots\,i_r}[P]]\,\varphi=
\sum^Q_{R=1}\msum{n}\Sb h_1,\,\ldots,\,h_\alpha\\k_1,\,\ldots,
\,k_\beta\endSb
\Theta^{\,h_1\ldots\,h_\alpha\, j_1\ldots\,j_s}_{k_1\ldots\,k_\beta\,i\,i_1\ldots\,i_r}[P,R]\ 
\vnabla^{k_1\ldots\,k_\beta}_{h_1\ldots\,h_\alpha}[R]\,\varphi.\quad
\mytag{16.3}
$$
And finally, in the case of an extended covector field $\bold X$ one
should write
$$
\hskip -2em
\gathered
[\nabla_i,\,\vnabla^{j_1\ldots\,j_s}_{i_1\ldots\,i_r}[P]]X_k=
-\sum^n_{h=1}D^{\,h\,j_1\ldots\,j_s}_{k\,i\,i_1\ldots\,i_r}[P]
\ X_h\ +\\
+\sum^Q_{R=1}\msum{n}\Sb h_1,\,\ldots,\,h_\alpha\\k_1,\,\ldots,
\,k_\beta\endSb
\Theta^{\ h_1\ldots\,h_\alpha\,j_1\ldots\,j_s}_{k_1\ldots
\,k_\beta\,i\,i_1\ldots\,i_r}[P,R]\ \vnabla^{k_1\ldots\,
k_\beta}_{h_1\ldots\,h_\alpha}[R]\,X_k.
\endgathered
\mytag{16.4}
$$
The components of $\bold D[P]$ and $\boldsymbol\Theta[P,R]$ in the
above three formulas \mythetag{16.2}, \mythetag{16.3}, \mythetag{16.4}
are taken from \mythetag{15.9} and \mythetag{15.10} respectively.\par
    The last commutation relationship is \mythetag{15.13}. In order to
write it in a local chart one should consider the commutator of two
covariant derivatives $\nabla_i$ and $\nabla_j$:
$$
\gather
\hskip -2em
\gathered
[\nabla_{\!i},\,\nabla_{\!j}]X^k=-\sum^n_{h=1}T^h_{ij}\ \nabla_{\!h}X^k
+\sum^n_{h=1}R^k_{h\,ij}\ X^h\ +\\
+\msum{n}\Sb h_1,\,\ldots,\,h_\alpha\\k_1,\,\ldots,\,k_\beta\endSb
\Omega^{\ h_1\ldots\,h_\alpha}_{\,ij\,k_1\ldots\,k_\beta }[R]
\ \vnabla^{k_1\ldots\,k_\beta}_{h_1\ldots\,h_\alpha}[R]\,X^k,
\endgathered
\mytag{16.5}\\
\vspace{2ex}
\hskip -2em
[\nabla_{\!i},\,\nabla_{\!j}]\varphi=-\sum^n_{h=1}T^h_{ij}\ 
\nabla_{\!h}\,\varphi+\msum{n}\Sb h_1,\,\ldots,\,h_\alpha\\
k_1,\,\ldots,\,k_\beta\endSb\Omega^{\ h_1\ldots\,h_\alpha}_{\,i
j\,k_1\ldots\,k_\beta }[R]\ \vnabla^{k_1\ldots\,k_\beta}_{h_1
\ldots\,h_\alpha}[R]\,\varphi,
\qquad
\mytag{16.6}\\
\vspace{2ex}
\hskip -2em
\gathered
[\nabla_{\!i},\,\nabla_{\!j}]X_k=-\sum^n_{h=1}T^h_{ij}\ \nabla_{\!h}X^k
-\sum^n_{h=1}R^h_{k\,ij}\ X_h\ +\\
+\msum{n}\Sb h_1,\,\ldots,\,h_\alpha\\k_1,\,\ldots,\,k_\beta\endSb
\Omega^{\ h_1\ldots\,h_\alpha}_{\,ij\,k_1\ldots\,k_\beta }[R]
\ \vnabla^{k_1\ldots\,k_\beta}_{h_1\ldots\,h_\alpha}[R]\,X_k,
\endgathered
\mytag{16.7}
\endgather
$$
The components of the torsion tensor $\bold T$, the components of the
curvature tensor $\bold R$, and the components of the tensor $\boldsymbol
\Omega[R]$ in \mythetag{16.5}, \mythetag{16.6}, \mythetag{16.7} are given
by the formulas \mythetag{15.17}, \mythetag{15.19}, and \mythetag{15.20}
respectively.\par
     The formulas \mythetag{16.2}, \mythetag{16.3}, \mythetag{16.4} and
\mythetag{16.5}, \mythetag{16.6}, \mythetag{16.7} are written for the 
cases of vectorial, covectorial, and scalar fields. However, the 
lemma~\mythelemma{7.5} and the theorem~\mythelemma{7.1} say that they are
sufficient for to write the analogous formulas in the case where the commutators $[\nabla_i,\,\vnabla^{j_1\ldots\,j_s}_{i_1\ldots\,i_r}[P]]$
and $[\nabla_{\!i},\,\nabla_{\!j}]$ are applied to the components 
of an arbitrary extended tensor field $\bold X$.\par
\head
17. Tensor functions of tensors\\
and the chain rule in tensorial form.
\endhead
    Tensor-valued functions with tensorial arguments appear rather often 
in applications. The most simple examples are the following ones:
\roster
\rosteritemwd=1pt
\item"---" the force field $\bold F(x^1,x^2,x^3,v^1,v^2,v^3)$ acting 
upon a point mass that moves according the Newton's second law;
\item"---" the Lagrange function $L(x^1,x^2,x^3,v^1,v^2,v^3)$ of such 
a point mass;
\item"---" the Hamilton function $H(x^1,x^2,x^3,p_1,p_2,p_3)$ of such
a point mass.
\endroster
These examples in a little bit more general form were mentioned in
section~2 (see comment to the formula \mythetag{2.2}) and in section~4.
Our next example is from the field theory. The action integral of the
electromagnetic field in vacuum is written as
$$
S=-\frac{1}{16\,\pi\,c}\int\sum^3_{i=0}\sum^3_{j=0}\sum^3_{\alpha=0}
\sum^3_{\beta=0}g_{ij}\ g_{\alpha\beta}\ F^{i\alpha}\,F^{j\beta}
\sqrt{-\det g\,}\,d^{\kern 0.5pt 4}\kern -0.5pt x
$$
(see \mycite{40} for details). The term under integration in this
formula is a scalar function
$$
\hskip -2em
L=-\frac{1}{16\,\pi\,c}\sum^3_{i=0}\sum^3_{j=0}\sum^3_{\alpha=0}
\sum^3_{\beta=0}g_{ij}\ g_{\alpha\beta}\ F^{i\alpha}\,F^{j\beta}.
\mytag{17.1}
$$
However, its value is determined by the tensor of the electromagnetic
field $\bold F$:
$$
F^{i\alpha}=\Vmatrix
0   &  -E^1 &  -E^2 &  -E^3\\
\vspace{1ex}
E^1 &   0   &  -H^3 & \,H^2\\
\vspace{1ex}
E^2 & \,H^3 &   0   &  -H^1\\
\vspace{1ex}
E^3 &  -H^2 & \,H^1 &   0\endVmatrix.
\mytag{17.2}
$$
Apart from \mythetag{17.2}, in \mythetag{17.1} we have the components
of the Minkowski metric:
$$
g_{ij}=\Vmatrix
1 &  0 &  0 &  0\\
0 & -1 &  0 &  0\\
0 &  0 & -1 &  0\\
0 &  0 &  0 & -1
\endVmatrix.
\mytag{17.3}
$$
In Cartesian coordinates the Minkowski metric is represented by the matrix
\mythetag{17.3}. If we use some curvilinear coordinate system, the matrix
components $g_{ij}$ become depending on the coordinates $x^0,\,x^1,\,x^2,\,
x^3$ of a point in the Minkowski space. In special relativity the role of
the Minkowski metric is not so significant as in general relativity. For
this reason, writing \mythetag{17.2} formally, we can indicate the presence
of $\bold g$ as an additional dependence on the spatial variables $x^0,\,
x^1,\,x^2,\,x^3$ in $L$:
$$
\hskip -2em
L=L(x^0,\ldots,x^3,F^{00},F^{01},F^{02},\ldots,F^{33}).
\mytag{17.4}
$$
For each particular configuration of the electromagnetic field $F^{ij}$ in
\mythetag{17.4} are some particular functions of $x^0,\,x^1,\,x^2,\,x^3$.
However, in some cases, e\.\,g\. in deriving the field equations from the
variational principle in form of the Euler-Lagrange equations, the quantities $F^{ij}$ are treated as independent variables.\par
     The example of the electromagnetic field, i\.\,e\. the function 
\mythetag{17.4}, can be considered as a background for various generalizations of the electromagnetism. Such theories could include
several tensorial fields $\bold T[1],\,\ldots,\,\bold T[Q]$. Therefore,
for the density in the action integral of such theories one should choose
some function $L$ depending on the variables \mythetag{3.3}:
$$
\hskip -2em
L=L(x^1,\ldots,x^n,\,T^{1\kern 1pt\ldots\,1}_{1\kern 1pt\ldots\,1}[1],
\ldots,T^{n\,\ldots\,n}_{n\,\ldots\,n}[Q]).
\mytag{17.5}
$$
This means that $L$ in \mythetag{17.5} is an extended scalar field
associated with some composite tensor bundle $N=T^{r_1\ldots\,r_Q}_{s_1\ldots\,s_Q}\!M$. If the whole scenario is
performed in the Minkowski space or in some space $M$ equipped with a
metric $\bold g$ and with some connection $\Gamma$, then the differentiations introduced in the
definition~\mythedefinition{12.2} and in the
definition~\mythedefinition{13.2} are applicable to $L$. On the
other hand, if some particular configuration of the fields 
$\bold T[1],\,\ldots,\,\bold T[Q]$ is given, then
$$
\cases T^{i_1\ldots\,i_r}_{j_1\ldots\,j_s}[1]=
T^{i_1\ldots\,i_r}_{j_1\ldots\,j_s}[1](x^1,\ldots,x^n)
\text{, \ where \ }r=r_1,\ s=s_1,\\
\ .\ .\ .\ .\ .\ .\ .\ .\ .\ .\ .\ .\ .\ .\ .
\ .\ .\ .\ .\ .\ .\ .\ .\ .\ .\ .\ .\ . 
\ .\ .\ .\ .\ .\ .\ .\ .\ .\ .\ .\ .\ .\\
T^{i_1\ldots\,i_r}_{j_1\ldots\,j_s}[Q]=
T^{i_1\ldots\,i_r}_{j_1\ldots\,j_s}[Q](x^1,\ldots,x^n)
\text{, \ where \ }r=r_Q,\ s=s_Q.
\endcases\quad
\mytag{17.6}
$$
Substituting \mythetag{17.6} into \mythetag{17.5}, we obtain
$$
\hskip -2em
\widetilde L=\widetilde L(x^1,\ldots,x^n).
\mytag{17.7}
$$
The function $\widetilde L$ in \mythetag{17.7} represents a standard (not
extended) scalar field. This means that we can differentiate $L$ in two
ways: as an extended field in its original form \mythetag{17.4} and as a standard field upon substituting some particular fields \mythetag{17.6}
into its arguments. The same is true for an arbitrary extended tensor field 
$\bold X$.
\mytheorem{17.1} Let\/ $\bold X$ be an extended tensor field of the 
type $(\alpha,\beta)$ associated with a composite tensor bundle
$N=T^{r_1\ldots\,r_Q}_{s_1\ldots\,s_Q}\!M$ and let $\bold T[1],\,\ldots,
\,\bold T[Q]$ be some non-extended tensor fields that determine some
particular section $q=q(p)$ of the bundle $N$. Denote by 
$\widetilde{\bold X}$ the non-extended tensor field obtained from $\bold X$
by substituting $\bold T[1],\,\ldots,\,\bold T[Q]$ into its arguments.
Then
$$
\hskip -2em
\nabla_{\bold Y}\widetilde{\bold X}=\nabla_{\bold Y}\bold X+
\sum^Q_{P=1}C(\nabla_{\bold Y}\bold T[P]\otimes
\vnabla[P]\,\bold X),
\mytag{17.8}
$$
where $\bold Y$ is some non-extended vector field in $M$, $\nabla_{\bold Y}\widetilde{\bold X}$ is the standard covariant
differentiation\footnotemark, $\nabla_{\bold Y}\bold X$ is the spacial
covariant differentiation, and $\nabla_{\bold Y}\bold T[P]$ is again the
standard covariant differentiation.
\endproclaim
\footnotetext{\ Writing covariant differentiations we assume that $M$ is 
equipped with some connection $\Gamma$. This can be either a standard connection or an extended connection. In the latter case it is converted 
to the standard connection by means of the section $q=q(p)$. In other 
words, one should substitute \mythetag{17.6} into the arguments of
$\Gamma^k_{ij}(x^1,\ldots,x^n,\,T^{1\kern 1pt\ldots\,1}_{1\kern 1pt\ldots\,1}[1],\ldots,T^{n\,\ldots\,n}_{n\,\ldots\,n}[Q])$.} 
\adjustfootnotemark{-1}
     The equality \mythetag{17.8} in the theorem~\mythetheorem{17.1} is a
tensorial form of the well-known chain rule for differentiating composite
functions. Its proof is pure calculations. First of all one should
write the equality \mythetag{17.8} in local coordinates. Here covariant
differentiations are replaced by covariant derivatives. As for the vector
field $\bold Y$, it can be dropped at all. As a result \mythetag{17.8} is
written as
$$
\hskip -2em
\gathered
\nabla_{\!i}\widetilde X^{i_1\ldots\,i_\alpha}_{j_1\ldots\,j_\beta}
=\nabla_{\!i}X^{i_1\ldots\,i_\alpha}_{j_1\ldots\,j_\beta}\,+\\
+\sum^Q_{P=1}
\msum{n}\Sb h_1,\,\ldots,\,h_r\\k_1,\,\ldots,\,k_s\endSb
\nabla_{\!i}T^{h_1\ldots\,h_r}_{k_1\ldots\,k_s}[P]\ \vnabla^{k_1\ldots
\,k_s}_{h_1\ldots\,h_r}[P]X^{i_1\ldots\,i_\alpha}_{j_1\ldots\,j_\beta}.
\endgathered
\mytag{17.9}
$$
The equality \mythetag{17.9} is derived by direct calculations based 
on the formulas \mythetag{13.15} and \mythetag{12.11}. The equality \mythetag{17.8} then is derived by multiplying both sides of 
\mythetag{17.9} by $Y^i$ and summing over the index $i$.

\enddocument
\end